\documentclass[11pt]{amsart}
\usepackage{latexsym}
\usepackage{amsmath,amssymb,amsthm}
\usepackage{amscd}
\usepackage[dvips]{graphics}
\newcommand{\Z}{\mathbb{Z}}
\newcommand{\R}{\mathbb{R}}

\newcommand{\cp}{\mathbb{CP}}
\newcommand{\zbar}{\overline{z}}
\newcommand{\del}{\partial}
\newcommand{\al}{\alpha}

\newcommand{\rp}{\mathbb{RP}}
\newcommand{\re}{\mathbb{R}}
\newcommand{\co}{\mathbb{C}}

\newcommand{\gl}[1]{\mathbf{GL}(#1,\mathbb{R})}
\newcommand{\na}{\nabla}
\newcommand{\sfrac}[2]{{\textstyle \frac{#1}{#2}}}

\newcommand{\D}{\displaystyle}

\newtheorem{prop}{Proposition} 

\newtheorem{thm}{Theorem}

\newtheorem{ex}{Example}
\newtheorem{lem}[prop]{Lemma}
\newtheorem{dfn}{Definition}

\theoremstyle{remark}
\newtheorem*{rem}{Remark}

\begin{document}
\title{Affine Manifolds, SYZ Geometry and The ``Y'' Vertex}
\author{John Loftin, Shing-Tung Yau, and Eric Zaslow}
\maketitle

\begin{abstract}
We prove the existence of a solution to the Monge-Amp\`ere
equation $\det {\rm Hess}(\phi)=1$ on a cone over a
thrice-punctured two-sphere. The total space of the tangent bundle
is thereby a Calabi-Yau manifold with flat special Lagrangian
fibers. (Each fiber can be quotiented to three-torus if the the
affine monodromy can be shown to lie in $\mathbf{SL}(3,\mathbb
Z)\ltimes \mathbb R^3$.) Our method is through Baues and
Cort\'es's result that a metric cone over an elliptic affine
sphere has a parabolic affine sphere structure (i.e., has a
Monge-Amp\`ere solution). The elliptic affine sphere structure is
determined by a semilinear PDE on $\mathbb{CP}^1$ minus three
points, and we prove existence of a solution using the direct
method in the calculus of variations.
\end{abstract}


\section{Introduction}

The basic question we would like to understand is, What does the
geometry of a Calabi-Yau manifold look like near (or ``at'') the
large complex structure limit point?  In order to answer this
question, one first fixes the ambiguity of rescaling the metric by
an overall constant.  Gromov proved that Ricci-flat manifolds with
fixed diameter have a limit under the Gromov-Hausdorff metric (on
the space of metric spaces).

Now by the conjecture of \cite{syz}, one expects that near the
limit, the Calabi-Yau has a fibration by special Lagrangian
submanifolds which are getting smaller and smaller (than the base).
The reason can be found by looking at the mirror large radius limit.
Fibers are mirror to the zero brane, and the base is mirror to the
$2n$ brane, which becomes large at large radius. Metrically, the
Calabi-Yau geometry should be roughly a fibration over the moduli
space of special Lagrangian tori ($T$). The dual fibration is by
dual tori, ${\rm Hom}\;(\pi_1(T),\R)/{\rm Hom}\;(\pi_1(T),\Z).$ The
flat fiber geometry of the dual torus fibration has a flat fiber
dual, which is not the same as the original geometry, but should be
the same after corrections by disk instantons.  These should get
small in the limit of small tori, though. Namely, we expect that the
Gromov-Hausdorff limit of a fixed-diameter Calabi-Yau manifold
approaching a maximal-degeneration point carries the same geometry
as the moduli space of special Lagrangian tori. That is, it is a
manifold (and an affine manifold at that) of half the dimension.

Further, the Calabi-Yau near the limit should be ``asymptotically
close'' to the standard flat torus fibration over special Lagrangian
tori moduli space, whose fibers are the flat tori (we get this from
dual torus considerations applied to the mirror manifold).  This
space is a quotient of the tangent space of the moduli space (which
is Hessian), and the Calabi-Yau condition means that the limiting
affine manifold metric should be Monge-Amp\`ere (${\rm det}\,{\rm
Hess}\,\Phi=1$). Global considerations require that the lattice
defining the torus (generated by the vectors associated to the
Hessian coordinates) is well-defined, meaning that the
Monge-Amp\`ere manifold has affine transition functions in the
semi-direct product of $\mathbf{SL}(n,\Z)$ with $\R^n$ translations.

This well-known conjecture (see e.g.\ \cite{GW} \cite{KS}
\cite{Fuk}) was proved by Gross-Wilson for the special case of K3
surfaces \cite{GW}.  Their proof uses the Ooguri-Vafa \cite{o-v}
metric in the neighborhood of a torus degeneration to build an
approximate Ricci-flat metric on the entirety of an elliptic K3
with 24 singular fibers (with an
elliptic-fibration/stringy-cosmic-string metric outside the
patches containing the degenerations).


Another aspect of this conjecture is that the limiting manifold
should have singularities in codimension two, with monodromy
transformations defined for each loop about the singular set.
Gross has shown the existence of a limiting singular set of
codimension two for (non-Lagrangian) torus fibrations on toric
three-fold Calabi-Yau's, and further has shown that the limiting
singular set has the structure of a trivalent graph \cite{gross}.
Also, W.D.\ Ruan has constructed Lagrangian torus fibrations with
codimension-two singular locus on quintic Calabi-Yau hypersurfaces
\cite{wd-ruan}. Taking our cue from these works, then in three
dimensions a point on the limiting manifold may be a smooth point,
a point near an interval singularity, or the trivalent vertex of a
``Y''-shaped singularity locus (these vertices have a
subclassification based on the monodromies near the vertex).
Examples of explicit Monge-Amp\`ere metrics for points of the
first two types are known, the interval singularity reducing to
the an interval times the two-fold point singularity. The absence
of a local metric model of a trivalent vertex singularity limits
our ability to prove this conjecture in three dimensions. Even if
we had such a model, it might not suffice to prove the conjectures
about the limiting metric, just as the non-Ooguri-Vafa elliptic
fibration metric does not suffice in the two-dimensional case.
Still, we regard the existence of a semi-flat Calabi-Yau metric
near the ``Y'' vertex as an important first step in addressing
these conjectures in three dimensions.\footnote{This vertex is not
the same as the topological vertex of \cite{AKMV}, which appears
at a corner of the toric polyhedron describing the Calabi-Yau. The
relation between the toric description of the Calabi-Yau and the
singularities of the special Lagrangian torus fibration has been
discussed in \cite{gs}.}

We should also remark that on other fronts, there has been much
progress recently in describing the proposed limit space.  There are
combinatorial constructions of integral affine manifolds with
singularities in the works of Haase-Zharkov \cite{hz02} and
Gross-Siebert \cite{gs}, who discuss mirror symmetry from
combinatorial and algebro-geometric points of view. Haase-Zharkov
\cite{hz03} also construct affine K\"ahler metrics on their
examples, but these do not satisfy the Monge-Amp\`ere equation.
Recently Zharkov has put forward a detailed conjectural picture of
the degeneration of Calabi-Yau metrics in toric hypersurfaces
\cite{zh}.

We therefore concern ourselves with studying Monge-Amp\`ere
manifolds in low dimensions, with the goal of finding a local model
for a trivalent degeneration of special Lagrangian tori. Taking our
model to be a metric cone over a thrice-punctured two-sphere, we
have, by an argument of Baues and Cort\'es, that the sphere metric
should be an elliptic affine sphere on $S^2$ with three
singularities.  The singularity type of the metric at the three
points is fixed by the pole behavior of a holomorphic cubic form.
Our main result is a proof of the existence of such an elliptic
affine sphere, hence its cone, which solves the Monge-Amp\`ere
equation with the desired singular locus.

The plan of attack is as follows. We study the Monge-Amp\`ere
equation and affine K\"ahler manifolds in Section 2, exhibiting a
few new solutions in three dimensions. In Section 3, we review
some basic notions in affine differential geometry, and recall
Baues and Cort\'es's result relating $n$-dimensional elliptic
affine spheres to $(n+1)$-dimensional parabolic affine spheres. In
Section 4, we study the relation between elliptic fibrations
(``stringy cosmic string'') and affine K\"ahler coordinates in two
dimensions. In Section 5, we recall Simon and Wang's theory of
two-dimensional affine spheres, focusing on the elliptic case.  In
Section 6, we first study the local structure near a singularity
of the elliptic affine sphere equation, and finally show the
existence of an elliptic affine sphere structure on $\cp^1$ minus
3 points.  This is our main result. The Monge-Amp\`ere metric near
the ``Y'' vertex is then constructed as a cone.

\begin{rem}
The key to finding an elliptic affine sphere metric on $\cp^1$ minus
3 points is the PDE
$$\psi_{z\bar z} + |U|^2e^{-2\psi} + \sfrac12 e^\psi=0,$$
where $U\,dz^3$ is a holomorphic cubic differential and
$e^\psi|dz|^2$ is the natural affine metric on an elliptic affine
sphere.  It is interesting to note that a similar equation,
$$\psi_{z\bar z} = e^{-2\psi} - e^\psi,$$
has also come up in the construction of special Lagrangian cones
in $\co^3$---see McIntosh \cite{mcintosh}.  In McIntosh's
construction, the equation comes from the geometry of
$\mathbf{SU}(3)$, while in the present case, the notion of an
elliptic affine sphere is invariant under $\mathbf{SL}(3,\re)$.
Such equations involving $e^{-2\psi}$ and $e^\psi$ go back to
\c{T}i\c{t}eica \cite{tzitz}, and are naturally associated to the
geometry of real forms of $\mathbf{SL}(3,\co)$.
\end{rem}

Correction: The first author would like to take this opportunity
to correct an erroneous attribution in \cite{loftin04}.  The
notion that parabolic affine spheres may be represented by
holomorphic data does not go back to Blaschke, but seems to be
originally due to Calabi \cite{calabi3}.  There is also a useful
Weierstrass type representation due to Ferrer-Mart\'inez-Mil\'an
\cite{ferrermm}.  The first author would like to thank Professor
Mart\'inez for pointing this out to him.

\section{Affine K\"ahler Metrics and the Monge-Amp\`ere Equation}

We recall that a metric is of \emph{Hessian} type if in
coordinates $\{x^i\}$ it has the form $ds^2 = \Phi_{ij}dx^i\otimes
dx^j,$ where $\Phi_{ij} =
\partial^2 \Phi/\partial x^i\partial x^j.$ Hitchin proved \cite{hit}
that natural metric (``McLean'' or ``Weil-Petersson'') on moduli
space of special Lagrangian submanifolds naturally has this
structure, and the semi-flat metric on the complexification (by flat
bundles) defined by the K\"ahler potential $\Phi$ is Ricci flat if
\begin{equation}\label{m-a}{\rm det}\;(\Phi_{ij}) = 1.\end{equation}

A manifold whose the coordinate gluing maps are all affine maps is
called an \emph{affine manifold}.  A Hessian metric on an affine
manifold is called \emph{affine K\"ahler}.  Note that our
definition of a Hessian metric is more general than that of an
affine K\"ahler metric; this distinction is not often made in the
literature.

In local Hessian coordinates we can compute the Christoffel
symbols $\Gamma^i{}_{jk} = \frac{1}{2}\Phi^{il}\Phi_{jkl}$ (where
$\nabla_j \del_k = \Gamma^i{}_{jk}\del_i$), and defining the
curvature tensor $R_{ij}{}^k{}_l = \partial_i
\Gamma^k{}_{jl}+\Gamma^k{}_{im}\Gamma^m{}_{jl} - (i\leftrightarrow
j)$ by  $[\nabla_i,\nabla_j]\del_k = R_{ij}{}^k{}_l \del_l,$ we
find
\begin{equation}
\label{riemcurv}
R_{ijkl} = -\frac{1}{4}\Phi^{ab}[\Phi_{ika}\Phi_{jlb}-
\Phi_{jka}\Phi_{ilb}].
\end{equation}


\subsection{Hessian Coordinate Transformations}

One asks, what coordinate transformations preserve the Hessian
form of the metric?  In particular, are there non-affine
coordinate changes which preserve the Hessian character of a given
metric?  If we try to write $ds^2 = \Phi_{ij}dx^i dx^j =
\Psi_{ab}dy^ady^b = \Psi_{ab}y^a{}_iy^b{}_j dx^idx^j,$ then the
consistency equations $\Phi_{ijk} = \Phi_{kji}$ yield conditions
on the coordinate transformation $y(x).$ Specifically, we have
$\partial_k (\Psi_{ab}y^a{}_i y^b{}_j ) =
\partial_i (\Psi_{ab}y^a{}_k y^b{}_j),$
which is equivalent to
$$\Psi_{ab}(y^a{}_{i}y^b_{jk} - y^a{}_{k}y^b{}_{ij}) = 0.$$

In two dimensions, for example, there can be many solutions to
these equations. In Euclidean space $\Psi_{ab} = \delta_{ab}$ with
coordinates $y^a,$ if we put $y^1 = f(x^1+x^2) + g(x^1-x^2)$ and
$y^2 = f(x^1 +x^2)-g(x^1-x^2),$ then the equations are solved and
we can find $\Phi(x).$ For example, if $f(s) = g(s) = s^2/2,$ we
find $\Phi(x) = [(x^1)^4 + 6x^1x^2 + (x^2)^4]/12.$

Note that this transformation is not affine.  Thus Hessian metrics
may exist on non-affine manifolds, and our notion of Hessian
metric is strictly broader than our notion of affine K\"ahler
metric. Affine K\"ahler manifolds can be characterized as locally
having an abelian Lie algebra of gradient vector fields acting
simply transitively \cite{ruuska}. Though Hessian manifolds are
not the same as affine manifolds, a Hessian manifold appearing as
a moduli space of special Lagrangian tori must have an affine
structure. We therefore focus on affine K\"ahler manifolds in this
paper.

\subsection{An Example of a Monge-Amp\`ere Metric}

As we will see in Section \ref{scs}, there are many Monge-Amp\`ere
metrics in two dimensions, but a paucity of examples in three or
more dimensions. Here we provide one detailed example and remark how
a few others may be found.

\begin{ex} \label{radial-sol}
In dimension $d$
consider the ansatz $\Phi = \Phi(r),$ where $r = \sqrt{\sum_i (x^i)^2}.$
As shown by Calabi, the equation (\ref{m-a}) is solved if
\begin{equation}\label{usol}
\Phi(r) = \int (1 + r^{d})^{1/d}.\end{equation}
The rescalings $r\rightarrow cr$ and $\Phi\rightarrow c\Phi$ also
have constant $\det({\rm Hess}\,\Phi).$
For example, in two dimensions ($d=2$),
$\Phi(r) = \sinh^{-1} (r) + r\sqrt{1+r^2}$
is a solution.
\end{ex}

\begin{rem}
It is also possible to find solutions in a few other cases in
dimension 3 by imposing symmetry.  In particular we let $\Phi$ take
the special forms for coordinates $x,y,z$ of $\re^3$.
\begin{itemize}
\item $\Phi = A(\rho)B(z)$ for $\rho=\sqrt{x^2+y^2}$.
\item $\Phi = \Phi(xyz)$.
\item $\Phi = \Phi(xy+yz+xz)$.
\end{itemize}
Solutions follow from straightforward ODE techniques.
\end{rem}

\section{Affine Spheres} Convex functions $\Phi$
satisfying the Monge-Amp\`ere equation $\det \Phi_{ij}=1$ have a
particularly useful interpretation in terms of affine differential
geometry.  The graph of such a $\Phi$ in $\re^{n+1}$ is a
\emph{parabolic affine sphere}. In this subsection, we introduce the
basic notions of affine differential geometry and recall a recent
result of Baues and Cort\'es which allows us to find 3-dimensional
solutions to the Monge-Amp\`ere equation by constructing a
2-dimensional elliptic affine sphere.

Affine differential geometry is the study of those properties of
hypersurfaces $H\subset \re^{n+1}$ which are invariant under
volume-preserving affine transformations. For basic background on
affine differential geometry, see Calabi \cite{calabi2}, Cheng-Yau
\cite{cheng-yau86} and Nomizu-Sasaki \cite{nomizu-sasaki}.  We
assume that $H$ is a smooth locally strictly convex hypersurface.
The affine normal $\xi$ to $H$ is a transverse vector field on $H$
which is invariant under the action volume-preserving affine
transformations on $\re^{n+1}$ in the sense that if
$\Psi\!:\re^{n+1}\to\re^{n+1}$ is such a transformation, then at all
$p\in H$,
$$\Psi_*(\xi_H(p)) = \xi_{\Psi(H)}(\Psi(p)).$$
We assume $\xi$ points inward (i.e., at $p$, $\xi(p)$ is on the same
side of any tangent plane of $H$ as $H$ is itself). Given such a
transverse vector field $\xi$, we have the following equations. Let
$X,Y$ be tangent vector fields on $H$ and let $D$ denote the
standard flat connection on $\re^{n+1}$.
\begin{eqnarray}
\label{dxy}
D_X Y &=& \na_X Y + h(X,Y)\xi, \\
\label{dxxi}
D_X \xi &=& -S(X).
\end{eqnarray}
Here $\na$ is a torsion-free connection on $H$, $h$ is a
Riemannian metric on $H$ (since $H$ is convex and $\xi$ points
inward), and $S$ is an endomorphism of the tangent space. $\na$ is
called the \emph{affine connection}, $h$ is the \emph{affine
metric}, and $S$ is the \emph{affine shape operator}.  The trace
of $S$ divided by the dimension $n$ is called the \emph{affine
mean curvature}.  Note that there is no part of equation
(\ref{dxxi}) in the span of $\xi$; if $D_X\xi\in T_p H$, then
$\xi$ is said to be \emph{equiaffine}.

The affine normal $\xi$ can be uniquely characterized as follows:
\begin{itemize}
\item $\xi$ points inward on $H$.
\item $\xi$ is equiaffine.
\item For any basis $X_1,\dots,X_n$ of the tangent space of $H$,
\begin{equation}\label{det-cond}
\det(X_1,\dots,X_n,\xi)^2 = \det h(X_i,X_j),
\end{equation}
where the determinant on the left is that on $\re^{n+1}$, and the
determinant on the right is that of an $n\times n$ matrix.
\end{itemize}

Another important invariant is the \emph{Pick form}, which may be
defined as the tensor which is the difference $C = \hat\na -  \na$
for $\hat\na$ the Levi-Civita connection of the affine metric.  The
Pick form satisfies the following \emph{apolarity} condition:
$$\sum_{i=1}^n C^i_{ij} = 0, \qquad j=1,\dots,n.$$ When the upper
index is lowered by the affine metric, the Pick form is totally
symmetric on all three indices.

A \emph{parabolic affine sphere} is hypersurface for which $\xi$
is a constant vector.  If $\xi=(0,\dots,0,1)$, then a parabolic
affine sphere can locally be written as a graph $(x,\Phi(x))$ for
a convex function $\Phi$ which satisfies the Monge-Amp\`ere
equation $\det \Phi_{ij}=1$. This condition may be checked by
using the conditions above for the affine normal.

A hypersurface $H$ is an \emph{elliptic affine sphere} if all the
affine normals point toward a given point in $\re^{n+1}$, called the
\emph{center} of $H$. In this case, the affine shape operator
$S=\lambda I$ for $\lambda>0$ and $I$ the identity operator on the
tangent space.  By translation, we may assume the center is the
origin, and by scaling, we may assume that $\lambda=1$.  In this
case, the affine normal $\xi$ is minus the position vector.

\begin{ex}
The unit sphere in $\re^{n+1}$ is an elliptic affine sphere centered
at the origin.  In this case, we may compute for $\xi$ equal to
minus the position vector that the affine metric $h$ is the
restriction of the Euclidean inner product.  It is straightforward
to check (\ref{det-cond}) is satisfied, and that $\xi$ is the affine
normal.
\end{ex}

Baues and Cort\'es establish a relationship between $n$-dimensional
elliptic affine spheres and $(n+1)$-dimensional parabolic affine
spheres \cite{BC}.  We use this result to reduce the problem of
finding 3-dimensional parabolic affine spheres to the problem of
finding 2-dimensional elliptic affine spheres.

\begin{thm}[Baues-Cort\'es] \label{bc-thm}
Let $H$ be an elliptic affine sphere in $\re^{n+1}$ centered at
the origin with affine mean curvature $1$.  Shrink $H$ if
necessary so that each ray through the origin hits $H$ only once.
Let $$\mathcal C = \bigcup_{r>0} rH$$ be the cone over $H$. Then
$\Phi = \frac12 {r^2}$ is convex and solves $\det \Phi_{ij}=1$ on
$\mathcal C$.  Using the diffeomorphism $\mathcal C \cong H \times
\re^+$, the affine K\"ahler metric satisfies
$$\frac{\partial^2 \Phi}{\partial x^i\partial x^j} \,dx^i\,dx^j =
r^2h + dr^2$$ for $h$ the affine metric on $H$.
\end{thm}

For the reader's convenience, we provide a proof of Baues and
Cort\'es's result, along the lines of \cite{loftin02}.

\begin{rem}
A similar proof shows that $$\Phi = \int (K r^{n+1} +
A)^{\frac1{n+1}}$$ solves $\det \Phi_{ij} = {\rm const.}$ for
constants $K$ and $A$.  In the case $H$ is the standard Euclidean
sphere in $\re^{n+1}$, we recover Calabi's Example \ref{radial-sol}
above.
\end{rem}

\begin{proof}
Assume $H$ is an elliptic affine sphere centered at 0 with affine
mean curvature 1, and form $\mathcal C$ and $\Phi$ from $H$ as
above. Denote the affine K\"ahler metric
$$g_{ij} = \frac{\partial^2 \Phi} {\partial x^i \partial x^j}.$$

Consider the position vector field $$X = x^i\frac{\partial}{\partial
x^i} = r\frac{\partial}{\partial r}.$$  Let $\Phi = r^2/2$. Note
that
$$X\Phi = x^i\frac{\partial \Phi}{\partial x^i} =
r\frac{\partial}{\partial r} \left(\frac12 r^2\right) = r^2 =
2\Phi.$$ Then take $\partial /\partial x^j$ to find
\begin{equation} \label{deriv}
x^i\,\frac{\partial^2\Phi}{\partial x^i \partial x^j} =
\frac{\partial \Phi}{\partial x^j}.
\end{equation}

Consider a vector $Y$ tangent to the hypersurface
$H=\left\{\Phi=\frac12\right\}$, so that $Y\Phi=0$. Then
(\ref{deriv}) shows that
\begin{equation} \label{xv}
 g(X,Y)=x^i\frac{\partial^2\Phi}{\partial
x^i \partial x^j} y^j = \frac{\partial \Phi}{\partial x^j} y^j =
Y\Phi =0,
\end{equation}
and also
$$g(X,X)=x^i\frac{\partial^2\Phi}{\partial x^i \partial x^j} x^j =
 X\Phi = r^2.$$

Now we'll show that the $g$ restricts to a multiple of the affine
metric on $H$.

Let $D$ be the canonical flat connection on $\re^{n+1}$.  Then our
affine K\"ahler metric $g$ is given by
\begin{equation} g(A,B)=(D_A d \Phi,B) \label{def-g} \end{equation}
where $A,B$ are vectors and $(\cdot\, ,\cdot)$ is the pairing
between one forms and vectors. $X$ is transverse to $H$. So at $x\in
H$, $\re^{n+1}=T_x(\re^{n+1})$ splits into $T_x(H) \oplus \langle X
\rangle$.  Then, since $-X$ is the affine normal,
\begin{equation} D_YZ= \na_YZ + h(Y,Z)(-X) \label{split}
\end{equation}
where $Y,Z$ are tangent vectors to $H$, $\na$ is a torsion-free
connection on $T(H)$, and $h$ is the affine metric.

Now consider
\begin{eqnarray*}
0&=&Y(d \Phi,Z)\\
&=&(d \Phi,D_YZ)+(D_Y d \Phi,Z)\\
&=&-r^2\,h(Y,Z) + g(Y,Z)
\end{eqnarray*}
by (\ref{xv}), (\ref{def-g}) and (\ref{split}). Therefore,
$g(Y,Z)=r^2\,h(Y,Z)$ for $Y,Z$ tangent to $H$.

So far, all calculations have been at a point in $H\subset \mathcal
C$.  A simple scaling argument shows:

\begin{prop} \label{perp}
Under the metric $g$, each level set of the potential $\Phi$ is
perpendicular to the radial direction $X$.  Under the diffeomorphism
$$\mathcal C \cong H \times \re^+,$$ we have
$$g(X,X)= r^2,\qquad g(Y,Z)=r^2 h(Y,Z),$$
where $Y$ and $Z$ are in the tangent space to $H$ and $h$ is the
affine metric of $H$.
\end{prop}

Now since $H$ is an elliptic affine sphere, we know that if
$Y_1,\dots,Y_n$ is a basis of the tangent space of $H$ at a point,
then
$$\det(Y_1,\dots,Y_n,-X)^2 = \det h(Y_i,Y_j).$$
Denote $-X$ by $Y_{n+1}$, and let $a,b$ be indices from $1$ to
$n+1$, while $i,j$ are indices from $1$ to $n$. Compute using
Proposition \ref{perp} for the standard frame on $\re^{n+1}$:
$$\det g_{ab} = \frac {\det g(Y_a,Y_b)}{\det(Y_1,\dots,Y_{n+1})^2} =
\frac {g(-X,-X) \cdot \det g(Y_i,Y_j)}{\det h(Y_i,Y_j)} = r^2\cdot
r^{2n} = r^{2n+2}.$$ Since on $H$, $r=1$, we have that for each
point on $H$, $$\det g_{ab} = \det \frac{\partial^2\Phi}{\partial
x^a\partial x^b} = 1.$$  For points not on $H$, the Monge-Amp\`ere
equation follows since $\Phi$ scales quadratically in $r$.
\end{proof}

\begin{ex}
If $H$ is the unit sphere in $\re^{n+1}$, then the corresponding
potential function $\Phi = \sfrac12 \|x\|^2$ clearly satisfies
$\det \Phi_{ij} = 1$ and the metric $\Phi_{ij}\,dx^i\,dx^j$ is the
standard flat metric on $\re^{n+1}$.
\end{ex}

\section{Two-Dimensional Monge-Amp\`ere Metrics and the
Stringy Cosmic String}
\label{scs}

Using a hyper-K\"ahler rotation we can treat any elliptic surface as
a special Lagrangian fibration and try to find its associated affine
coordinates and---if the fibers are flat---the corresponding
solution to the Monge-Amp\`ere equation.

In the case of the stringy
cosmic string, we begin with a semi-flat fibration with
torus fiber coordinates $t \sim t + 1$ and
$x\sim x+1.$  As a holomorphic
fibration, the stringy cosmic string is defined by a holomorphic
modulus $\tau(z).$  One can derive the K\"ahler potential
through the Gibbons-Hawking ansatz (with $\partial/\partial t$
as Killing vector) using connection one-form $A = -\tau_1 dx$
and potential $V = \tau_2$ (so $*dA = dV$),
then solving for the holomorphic coordinate.
One finds $\xi = t + \tau(z)x =
t + \tau_1 x + i \tau_2 x.$
The hyper-K\"ahler structure is specified by the forms
\begin{eqnarray*}
\omega_1 &=& dt \wedge dx +
(\frac{i}{2})\tau_2 dz\wedge d\overline{z}\\
\omega_2 + i\omega_3 &=& dz \wedge d\xi.
\end{eqnarray*}
The stringy cosmic string solution
starts directly from the K\"ahler potential
$K(z,\xi) = \xi_2/\tau_2 + k(z,\overline{z}),$ where $\partial_z
\partial_{\overline{z}}k = \tau_2.$

We seek the affine coordinates for the base of the semi-flat special
Lagranian torus fibration.  In coordinates $(x,t,z_1,z_2)$ the
metric has the block diagonal form
\begin{equation}\label{metblock}
Q\oplus R \equiv
\frac{1}{\tau_2}\left(\begin{array}{cc}|\tau|^2&\tau_1\\ \tau_1& 1
\end{array}\right)\oplus
\left(\begin{array}{cc}\tau_2&0\\0&\tau_2\end{array}\right).
\end{equation}
For a semi-flat fibration over an affine K\"ahler manifold in affine
coordinates, the base-dependent metric on the fiber looks the same
as the metric on the base. Therefore, we would need to find
coordinates $u_1(z,\zbar), u_2(z,\zbar)$ so that the metric in
$u$-space looks like $Q$ in (\ref{metblock}). This is accomplished
if the change of basis matrix $M_{ij} = \del z_i/\del u_j$ obeys
$M^T M = Q/\tau_2.$ A calculation reveals the general solution to be
$M = O\widetilde{M},$ where $\widetilde{M} =
\frac{1}{\tau_2}\left(\begin{array}{cc} \tau_2&0\\
\tau_1&1\end{array}\right),$ and $O$ is an orthogonal matrix. The
same result can be obtained using Hitchin's method, which we now
review.

Hitchin \cite{hit} obtains affine coordinates on the moduli space of
special Lagrangian submanifolds from period integrals.  In order to
apply this technique here we first make a hyper-K\"ahler rotation,
so that the fibration is special Lagrangian, by putting $\omega =
\omega_2,$ ${\rm Im}\,\Omega = \omega_3.$ ($\omega$ and $\Omega$ are
the symplectic  and holomorphic forms of the Calabi-Yau metric,
respectively.) Explicitly, for each base coordinate $z_i$ we
construct closed one forms on the Lagrangian $L,$ defined by
$\iota_{\del/\del z_i}\omega = \theta_i$ and compute the periods
$\lambda_{ij} = \int_{A_i}\theta_j,$ where $\{A_i\}$ is a basis for
$H_1(L,{\mathbb Z}).$ In our case we readily find $\theta_1 = dt +
\tau_1 dx,$ $\theta_2 = -\tau_2 dx,$ and, using the basis $A_1 = \{t
\rightarrow t + 1\},$ $A_2 = -\{x\rightarrow x+1\},$ we get
$$\lambda_{ij} = \left(\begin{array}{cc}1&0\\
-\tau_1&\tau_2\end{array}\right).$$
The forms $\lambda_{ij}dz_j$ are closed on
the base and we set them equal to $du_i.$
This defines the coordinates $du_i$ up to constants, and we find
$u_1 = z_1,$ $u_2 = -{\rm Re}\phi,$
where $\del_z \phi = \tau.$
To connect with the solution above, one easily
inverts the matrix  $\lambda_{ij} = \del u_i/\del z_j$
to find the matrix $\del z_i/\del u_j = \widetilde{M}.$


Legendre dual coordinates $v_i$ are defined as follows. Define
$(d-1)$-forms $\psi_i$ (here $d=2$ so the $\psi$'s are also
one-forms) by putting $\iota_{\del/\del z_i}{\rm Im}\,\Omega =
\psi_i$ and compute the periods $\mu_{ij} = \int_{B_i}\psi_j,$ where
$B_j \in H_{d-1}(L,{\mathbb Z})$ are Poincar\'e dual to the $A_i.$
In our example, $\psi_1 = \tau_2 dx,$ $\psi_2 = dt + \tau_1 dx,$
$B_1 = \{ x\rightarrow x+1\},$ $B_2 = \{t\rightarrow t + 1\},$ and
$$\mu_{ij} = \left(\begin{array}{cc}\tau_2 &\tau_1 \\
0&1\end{array}\right).$$
(Note $\lambda^T\mu$ is symmetric, as required.)
Setting $dv_i = \mu_{ij} dz_j$ we find
$v_1 = {\rm Im}\phi$ and $v_2 = z_2.$

Hitchin showed that the coordinates $u_i$ and $v_i$ are related by
the Legendre transformation defined by the function $\Phi$ whose
Hessian gives the metric.  Namely, $v_i = \del\Phi/\del u_i.$  We
can think of $\Phi$ as a function of the $z_i(u_j)$ and
differentiate with respect to $u_j$ using the chain rule. (We find
$\del z_i/\del u_j$ by inverting the matrix of derivatives $\del
u_i/\del z_j.$) One finds
$$\Phi_{z_1} = {\rm Im}\; \phi - \tau_1 z_2, \qquad \Phi_{z_2}
= \tau_2 z_2.$$ (For the transformed potential $\Psi$ we have
$\Psi_{z_1} = \tau_2 z_1$ and $\Psi_{z_2} = \tau_1 z_1 - {\rm
Re}\phi.$) The solution can be given in terms of another
holomorphic antiderivative,\footnote{V. Cort\'es has found a
generalization of this potential as the defining function of
special K\"ahler manifolds, which are locally special examples of
parabolic affine spheres in even dimensions, described by
holomorphic data \cite{cortes}. In the present case of dimension
2, Monge-Amp\`ere metrics were described using holomorphic data by
Calabi \cite{calabi3} and Ferrer-Mart\'inez-Mil\'an
\cite{ferrermm}.  Of course all these descriptions of
two-dimensional parabolic affine spheres using holomorphic data
are equivalent.} $\chi,$ such that $\del_z \chi = \phi.$
$$\Phi\; =\; - z_2 {\rm Re} \phi \;+\; {\rm Im}\chi.$$
Note, then, that being able to write down the explicit affine
K\"ahler potential depends only on our ability to integrate $\tau$
and invert the functions $u_i(z_j).$  The Legendre-transformed
potential is $\Psi = z_1{\rm Im}\phi - {\rm Im}\chi.$

\begin{ex} $\tau = 1/z.$
If we put
$z = re^{i\theta}$ and take $\tau = 1/z$ then $\phi = \log z,$ so
$u_1 = z_1$ and $u_2 = -{\rm Re}\phi =
-\log r.$  Thus $\Phi(u_1,u_2) = - z_2 \log \sqrt{z_1^2 + z_2^2} +
\int \log\sqrt{z_1 ^2 + z_2^2}dz_2,$ where
$z_1 = u_1$ and $z_2 = \sqrt{e^{-2u_2} - u_1^2}.$
Since $v_1 = {\rm Im}\Phi = \tan^{-1}(z_2/z_1)$ and $v_2 = z_2,$
we may solve the equations $\del\Phi/\del u_i = v_i$
to find
$$\Phi = u_1\left[\tan^{-1}\left(\sqrt{(e^{-u_2}/u_1)^2-1}\right) -
\sqrt{(e^{-u_2}/u_1)^2-1}\right].$$
One easily checks that $\det(\Phi_{ij})=1.$
\end{ex}


To summarize, let $z$ be a holomorphic coordinate on the base of a
semi-flat elliptic fibration. Let $\tau=\tau(z)$ be the
holomorphically varying modulus of the elliptic curve on the
fiber. Then we define $\phi,\chi$ holomorphic so that
$$ \phi_z=\tau, \qquad \chi_z=\phi.$$
Let $z=z_1+iz_2$ represent real and imaginary parts,
with similar notation for the real and imaginary parts of $\tau,\phi,\chi$.
Then affine flat coordinates $u_1,u_2$ may be chosen as
$$u_1=z_1, \qquad u_2= -\phi_1.$$
The metric on the base is given by
$$ \tau_2|dz|^2 =
\frac{\partial^2\Phi}{\partial u_i \partial u_j}\,du_idu_j$$
for the affine K\"ahler potential $\Phi$, which satisfies
$$\Phi = -z_2 \phi_1 + \chi_2.$$
The Legendre dual coordinates $v_i=\partial\Phi / \partial u_i$ are given by
$$v_1=\phi_2,\qquad v_2=z_2.$$
The potential $\Psi$ in the $v$ coordinates is the Legendre transform
of $\Phi$:
$$\Psi= u_1v_1+u_2v_2 - \Phi= z_1\phi_2-\chi_2.$$
$\Phi$ and $\Psi$ satisfy the Monge-Amp\`ere equation
$$\det \left(\frac{\partial^2\Phi}{\partial u_i \partial u_j}\right) = 1,
\qquad
\det \left(\frac{\partial^2\Psi}{\partial v_i \partial v_j}\right)=1.$$
The metric satisfies
$$ \tau_2|dz|^2 = \frac{\partial^2\Phi}{\partial u_i \partial u_j}\,du_idu_j
=\frac{\partial^2\Psi}{\partial v_i \partial v_j}\,dv_idv_j.$$

\section{Simon and Wang's Developing Map}
\label{devmap}

U.\ Simon and C.P.\ Wang \cite{simon-wang} formulate the condition
for a two-dimensional surface to be an affine sphere in terms of
the conformal geometry given by the affine metric. Since we rely
heavily on this work, we give a version of the arguments here for
the reader's convenience.  We are primarily interested in
constructing three-dimensional parabolic affine spheres by writing
them as cones over elliptic affine spheres in dimension two by
using Baues and Cort\'es's  Theorem \ref{bc-thm}.  Therefore, we
focus our discussion to the case of elliptic affine spheres in
dimension two, and conclude with some remarks about
two-dimensional parabolic affine spheres from this point of view.
(In Section \ref{solve-eas}, it will be useful to compare elliptic
and parabolic affine spheres in dimension two, particularly since
two-dimensional parabolic affine spheres admit exact solutions.)

\subsection{Elliptic Affine Spheres}

Before we get into the construction, a few remarks are in order.  We
consider a parametrization $f\!: \mathcal D\to \re^3$ where
$\mathcal D \subset \co$ is simply connected and $f$ is conformal
with respect to the affine metric.  Simon and Wang's procedure
involves writing the structure equations of the affine sphere as a
first-order system of PDEs (an initial-value problem) in the frame
$\{f,f_z,f_{\bar z}\}$---equations (\ref{z-deriv}-\ref{zbar-deriv})
below. By the Frobenius Theorem, this initial-value problem can be
solved as long as certain integrability conditions are satisfied.
One of these integrability conditions is a semilinear elliptic PDE
in the conformal factor of the affine metric---equation
(\ref{psi-eq}) below. Solving this PDE on a Riemann surface $\Sigma$
then provides an immersion from the universal cover $\tilde \Sigma$
to $\re^3$, the image being an (immersed) elliptic affine sphere. We
call this immersion Simon and Wang's \emph{developing map}.

In the the case of elliptic affine spheres, we take $\Sigma =
\cp^1$ minus 3 points. Integrating the initial value problem along
a path in $\pi_1\Sigma$ computes the monodromy of an
$\rp^2$-structure on $\Sigma$, which upon applying Baues and
Cort\'es's cone construction, provides the monodromy of the affine
flat structure on the cone over $\Sigma$, which is $\re^3$ minus a
``Y" vertex topologically.

We have not yet completed the ODE computation of the monodromy in
the present case of an elliptic affine sphere, but note that this
approach has been used to compute monodromy for convex $\rp^2$
structures (using hyperbolic affine spheres) \cite{loftin03}, and
also for a global existence result for parabolic affine spheres on
$S^2$ minus singular points \cite{loftin04}.

Consider a 2-dimensional elliptic affine sphere in $\re^3$.  Then
the affine metric gives a conformal structure, and we choose a local
conformal coordinate $z=x+iy$ on the hypersurface. The affine metric
is given by $h=e^{\psi}|dz|^2$ for some function $\psi$. Parametrize
the surface by $f:\mathcal{D} \rightarrow \mathbb{R}^3$, with
$\mathcal{D}$ a domain in $\mathbb{C}$. Since
$\{e^{-\frac{1}{2}\psi} f_x,e^{-\frac{1}{2}\psi} f_y\}$ is an
orthonormal basis for the tangent space, the affine normal $\xi$
must satisfy this volume condition (\ref{det-cond})
\begin{equation} \label{detxy}
\det(e^{-\frac{1}{2}\psi} f_x,e^{-\frac{1}{2}\psi} f_y, \xi) = 1,
\end{equation}
 which implies
\begin{equation}
\det(f_z,f_{\bar{z}},\xi)=\sfrac{1}{2}i e^\psi. \label{det-eq}
\end{equation}

Now only consider elliptic affine spheres centered at the origin and
with affine mean curvature scaled to be 1. In this case, the affine
normal is $-f$ (minus the position vector) and we have
\begin{equation}
\left\{ \begin{array}{c}
D_X Y = \na_X Y + h(X,Y)(-f)\\
D_X (-f) = -X
\end{array} \right. \label{struc}
\end{equation}
Here $D$ is the canonical flat connection on $\re^3$, $\na$ is a
torsion-free connection on the affine sphere, and $h$ is the affine
metric.

It is convenient to work with complexified tangent vectors, and we
extend $\na$, $h$ and $D$ by complex linearity.  Consider the
frame for the tangent bundle to the surface $\{ e_1 = f_z =
f_*(\frac{\partial}{\partial z}), e_{\bar 1}=f_{\bar z} =
f_*(\frac{\partial}{\partial {\bar z}}) \}$. Then we have
\begin{equation}
h(f_z,f_z)=h(f_{\bar z}, f_{\bar z})=0, \quad h(f_z,f_{\bar z}) =
\sfrac{1}{2}e^\psi.
\label{h-met}
\end{equation}
Consider $\theta$ the matrix of connection one-forms
$$\na e_i = \theta^j_i e_j, \quad i,j \in \{1,{\bar 1}\},$$
and $ {\hat \theta} $ the matrix of connection one-forms for the
Levi-Civita connection.   By (\ref{h-met})
\begin{equation}
{\hat \theta}^1_{\bar 1} ={\hat \theta}^{\bar 1}_1 = 0, \quad
{\hat \theta}^1_1 = \partial \psi, \quad
{\hat \theta}^{\bar 1}_{\bar 1} = {\bar \partial} \psi. \label{levi-cit}
\end{equation}

The difference ${\hat \theta} - \theta$ is given by the Pick form.  We
have
$${\hat \theta}^j_i - \theta^j_i = C^j_{ik} \rho^k,$$
where $\{ \rho^1 = dz,\rho^{\bar 1} = d{\bar z} \}$ is the dual
frame of one-forms.   Now we differentiate (\ref{det-eq}) and use
the structure equations (\ref{struc}) to conclude
$$\theta^1_1 + \theta^{\bar 1}_{\bar 1} = d \psi.$$
This implies, together with (\ref{levi-cit}), the apolarity condition
$$C^1_{1k} + C^{\bar 1}_{{\bar 1} k} = 0, \quad k \in \{1,{\bar 1} \}.$$
Then, when we lower the indices, the expression for the metric
(\ref{h-met}) implies that
$$C_{{\bar 1}1k} + C_{1{\bar 1}k} = 0.$$
Now $C_{ijk}$ is totally symmetric on three indices
\cite{cheng-yau86,nomizu-sasaki}. Therefore, the previous equation
implies that all the components of $C$ must vanish except
$C_{111}$ and $C_{ {\bar 1}{\bar 1}{\bar 1}} =
\overline{C_{111}}$.

This discussion completely determines $\theta$:
\begin{equation}
\left( \begin{array}{cc}  \theta^1_1 & \theta^1_{\bar 1} \\[1mm]
                                \theta^{\bar 1}_1 & \theta^{\bar 1}_{\bar1}
                  \end{array} \right)
 = \left( \begin{array}{cc}  \partial \psi & C^1_{{\bar 1}{\bar 1}}
                                d{\bar z} \\[1mm]
                                C^{\bar 1}_{11} dz & \bar{\partial} \psi
                  \end{array} \right)
=\left( \begin{array}{cc}  \partial \psi & \bar{U} e^{-\psi} d{\bar z} \\
                        U e^{-\psi} dz & \bar{\partial} \psi
                  \end{array} \right),
\label{conn-eq}
\end{equation}
where we define $U = C^{\bar 1}_{11} e^\psi$.

Recall that  $D$
is the canonical flat connection induced from ${\mathbb R}^3$.  (Thus,
for example, $D_{f_z}f_z = D_{\frac{\partial}{\partial z}} f_z =
f_{zz}$.)
Using this statement, together with (\ref{h-met}) and (\ref{conn-eq}),
the structure equations (\ref{struc}) become
\begin{equation}
\left\{ \begin{array}{c@{\,\,\,\,=\,\,\,\,}l}
f_{zz} & \psi_z f_z + U e^{-\psi} f_{\bar z} \\
f_{{\bar z}{\bar z}} & {\bar U} e^{-\psi} f_z + \psi_{\bar z}
f_{\bar z}
\\
f_{z{\bar z}} & -\frac{1}{2}e^\psi f \end{array} \right.
\label{fzz-eq}
\end{equation}
Then, together with the equations $(f)_z=f_z$, $(f)_{\bar z} =
f_{\bar z}$, these form a linear first-order system of PDEs in $f$,
$f_z$ and $f_{\bar z}$:
\begin{eqnarray} \label{z-deriv}
\frac{\partial}{\partial z} \left( \begin{array}{c}
f \\ f_z \\
f_{\bar z} \end{array} \right) &=& \left( \begin{array}{ccc}
0&1&0 \\
0& \psi_z & U e^{-\psi} \\
-\frac12 e^\psi & 0 & 0
\end{array} \right)
\left( \begin{array}{c} f \\ f_z \\ f_{\bar z} \end{array} \right),
\\ \label{zbar-deriv}
\frac{\partial}{\partial \bar z} \left( \begin{array}{c} f
\\f_z
\\ f_{\bar z} \end{array} \right) &=& \left( \begin{array}{ccc}
0&0&1 \\
-\frac12 e^\psi & 0 & 0 \\
0 & \bar{U} e^{-\psi} & \psi_{\bar z} \\
\end{array} \right)
\left( \begin{array}{c} f \\ f_z \\ f_{\bar z} \end{array} \right).
\end{eqnarray}
In order to have a solution of the system (\ref{fzz-eq}),
the only condition is
that the mixed partials must commute (by the Frobenius theorem). Thus we
require
\begin{eqnarray}
\psi_{z {\bar z}} + |U|^2 e^{-2\psi} +\sfrac12 e^\psi &=& 0,
\label{psi-eq} \\
\nonumber
U_{\bar z} &=& 0.
\end{eqnarray}

The system (\ref{fzz-eq}) is an initial-value problem, in that given
(A) a base point $z_0$, (B) initial values $f(z_0)\in\re^3$,
$f_z(z_0)$ and $f_{\bar z}(z_0)=\overline{f_z(z_0)}$, and (C) $U$
holomorphic and $\psi$ which satisfy (\ref{psi-eq}), we have a
unique solution $f$ of (\ref{fzz-eq}) as long as the domain of
definition $\mathcal{D}$ is simply connected.  We then have that the
immersion $f$ satisfies the structure equations (\ref{struc}). In
order for $-f$ to be the affine normal of $f(\mathcal{D})$, we must
also have the volume condition (\ref{det-eq}), i.e.\ $\det (f_z,
f_{\bar z}, -f)=\frac{1}{2}ie^\psi$.  We require this at the base
point $z_0$ of course:
\begin{equation}
\det(f_z(z_0), f_{\bar z}(z_0), -f(z_0)) =
\sfrac{1}{2}ie^{\psi(z_0)}. \label{init-val}
\end{equation}
Then use (\ref{fzz-eq}) to show that the derivatives with respect to
$z$ and ${\bar z}$ of $\det(f_z,f_{\bar z},-f)e^{-\psi}$ must
vanish. Therefore the volume condition is satisfied everywhere, and
$f(\mathcal{D})$ is an elliptic affine sphere with affine normal
$-f$.

Using (\ref{fzz-eq}), we compute
\begin{equation} \label{U-f}
\det(f_z,f_{zz},-f)=\sfrac12iU,
\end{equation}
which implies that $U$ transforms as a section of $K^3$, and
$U_{\bar z}=0$ means it is holomorphic.

Note that equation (\ref{psi-eq}) is in local coordinates. In
other words, if we choose a local conformal coordinate $z$, then
the Pick form ${\bf U}=U\,dz^3$, and the metric is
$h=e^\psi|dz|^2$. Then plug $U,$ $\psi$ into (\ref{psi-eq}).  In a
patch with a new holomorphic coordinate $w(z),$ the metric will
have the form $e^{\widetilde{\psi}}|dw|^2,$ with cubic form
$\widetilde{U}dw^3.$ Then $\widetilde{\psi}(w), \widetilde{U}(w)$
will satisfy (\ref{psi-eq}).

\subsection{Parabolic Affine Spheres}
Here we very briefly recall analogues of some results of
the last subsection for two-dimensional parabolic affine spheres.
This is due to Simon-Wang \cite{simon-wang}, and there is a
derivation similar to the one above in \cite{loftin04}.

A smooth, strictly convex hypersurface $H$ is a parabolic affine
sphere if the affine normal $\xi$ is a constant vector.  In $\re^3$,
we let $\xi=(0,0,1)$.  Let $f\!:\mathcal D\to \re^3$ be an immersion
of the parabolic affine sphere which is conformal with respect to
the affine metric.  Then $\{\xi,f_z,f_{\bar z}\}$ is a complexified
frame of $\re^3$ at each point in $H$.  The affine structure
equations lead to an initial-value problem similar to equations
(\ref{z-deriv}-\ref{zbar-deriv}) above, and the integrability
conditions are $U_{\bar z}=0$ (for the Pick form $U$) and
\begin{equation}
\label{inteq} \psi_{z\bar z} + |U|^2 e^{-2\psi} = 0
\end{equation}
for the affine metric $e^\psi|dz|^2$. This equation has many easy
explicit solutions, which is an advantage over the corresponding
equation (\ref{psi-eq}) for elliptic affine spheres.  In particular,
we often treat the extra term $\frac12 e^\psi$ in (\ref{psi-eq}) as
a perturbation of equation (\ref{inteq}).

It is also possible to find explicit solutions to the initial
value problem for two-dimensional parabolic affine spheres using
ODE techniques. This is unsurprising, as we saw in Section
\ref{scs} that the structure equations are completely integrable.
For example, if $U=z^{2\alpha-3}\,dz^3$, $z=\rho e^{i\theta}$, and
$A,B,\xi\in\re^3$ satisfy $\det(A,B,\xi)=-8$, then the following
$f$ is an immersed parabolic affine sphere in $\re^3$ with affine
normal $\xi$:
\begin{eqnarray*}
f &=& \sfrac 1{2\alpha} A \rho^\alpha \left[ \theta \cos
\alpha\theta -
\sfrac1\alpha \sin \alpha\theta + (\log \rho)\sin\alpha x \right] \\
&&{} -\sfrac1{2\alpha} B \rho^\alpha \cos \alpha \theta +
\sfrac1{\alpha^2} \xi \rho^{2\alpha} \left[ \sfrac 1{2\alpha} \cos
2\alpha\theta + \sfrac1\alpha - \log \rho \right].
\end{eqnarray*}
Note this solution has nontrivial monodromy around $z=0$.

\begin{rem}
In equation (\ref{inteq}), the transformation $\varphi = \log |U| -
\frac{1}{2}\psi$ results in the condition that $e^{2\varphi}|dz|^2$
has constant curvature $-4.$ We thank R.\ Bryant for pointing this
out.
\end{rem}

\section{Elliptic Affine Spheres and the ``Y'' Vertex}
\label{solve-eas}

In this section we will prove the existence of elliptic affine
two-sphere metrics with singularities---first locally near a
singularity (we find a radially symmetric solution), and then
globally on $S^2$ minus three points. The metric cone yields a
parabolic affine sphere metric near the ``Y'' vertex.

\subsection{Local Analysis}
\label{local}

Recall that given a holomorphic cubic differential $U$ on a domain
in $\co$ with coordinate $z$, a solution $\psi$ to
$$ \psi_{z\zbar}
+ |U|^2e^{-2\psi} + \sfrac{1}{2} e^\psi = 0
$$
provides an affine metric $e^\psi|dz|^2$ for an elliptic affine
sphere. The elliptic affine sphere can be reconstructed by Simon and
Wang's developing map.

Since Baues and Cort\'es result gives a parabolic affine sphere on
the cone over an elliptic sphere, a solution to this equation on the
thrice-punctured sphere will lead to a parabolic affine sphere on
$\R^3$ minus a Y-shaped set---whence a semiflat special Lagrangian
torus fibration over this base.  In the present subsection, we prove
the existence of radially symmetric solutions to (\ref{psi-eq}),
while we discuss the more global setting of the thrice-punctured
sphere in the next subsection.

For definiteness, we consider the case $U = z^{-2}$ and make the
ansatz $\psi = \psi(|z|).$  We look near $z=0,$ so we make the
change of variables $t = -\log |z|,$ $t\in (T,\infty),$ $T\gg0.$
This leads to the equation
\begin{equation}
\label{ellaffpsi}
N(\psi) := \partial_t^2 \psi + 4e^{-2(\psi-t)} + 2e^{\psi -2t} =0.
\end{equation}
We put $\psi = \psi_0 + \phi,$ where $\psi_0 = t + \log(2t)$ is the
solution to the parabolic equation (\ref{inteq}). Note that the last
term in (\ref{ellaffpsi}) is $O(te^{-t})$ for this function. We want
to solve $N(\psi_0+\phi)=0,$ which we expand as
\begin{equation}
\label{maineq}
N(\psi_0+\phi) =
N(\psi_0) + dN(\phi)\vert_{\psi_0} + Q(\phi)\vert_{\psi_0},
\end{equation}
where $Q(\phi)$ contains quadratic and higher terms.
Explicitly,
$$Q(\phi) = \frac{1}{t^2}\left(e^{-2\phi}-(1-2\phi)\right)
+4te^{-t}\left(e^\phi-(1+\phi)\right).$$
Note that $Q$ is not even a differential operator.
One calculates
$$N(\psi_0) = 4te^{-t},$$
$$dN(\phi)\vert_{\psi_0}
=: L\phi := \left[\partial_t^2 + V(\psi_0)\right]\phi,$$
where
$$V(\psi_0) = -8e^{-2(\psi_0 - t)} + 2e^{\psi_0-2t} = -\frac{2}{t^2}+
4te^{-t}.$$ Thus $L\phi = (\partial_t^2 - \frac{2}{t^2} + 4te^{-t})
\phi = (L_0 + 4te^{-t})\phi,$ where $L_0 = \partial_t^2 -
\frac{2}{t^2}.$ The equation (\ref{psi-eq}) is now
$$L\phi = f - Q(\phi),$$
with $f = -4te^{-t}.$ The idea will be to find an appropriate
Green function $G$ for $L,$ in terms of which a solution to this
equation becomes a fixed point of the mapping $\phi \rightarrow
G(f - Q(\phi))$---then to find a range of $\phi$ where this is a
contraction map, whence a solution by the fixed point theorem.

We claim that this map is a contraction for $\phi \sim O(te^{-t}).$
More specifically, consider for a value of $T>2$ to be determined
later, the Banach space $\mathcal B$ of continuous functions on
$[T,\infty)$ with norm  $$\|g\|_{\mathcal B} = \sup_{t\ge T}
\frac{g(t)}{te^{-t}}.$$ Showing  the map $\phi\rightarrow
G(f-Q(\phi))$ is a contraction map  involves estimating $Gf$ and
$GQ\phi.$ In fact, since $Q$ is a quadratic, non-derivative
operator, it is easy to see that $Q\phi$ is order $te^{-t}$ (even
smaller).  We then show that $G$ preserves the condition
$O(te^{-t})$ by showing $Gf\in\mathcal B$ (recall that $f\in\mathcal
B$, too). To find $G,$ we write $L = L_0 - \delta_L,$ where
$\delta_L = -4te^{-t},$ so that $G = L^{-1} = L_0^{-1} +
L_0^{-1}\delta_L L_0^{-1} + ... .$  To solve the equation $Lu = f,$
we first note that the change of variables $v = u+1$ leads to the
equation $Lv = -\frac{2}{t^2}.$ Let $v_0 =
L_0^{-1}(-\frac{2}{t^2})=1.$ Then define $v_{k+1} = L_0^{-1}
\delta_L v_k.$  Then $v = \sum_{k=0}^{\infty} v_k$ and $Gf = u =
\sum_{k=1}^{\infty} v_k.$

\begin{lem}
$|v_k(t)| < (16te^{-t})^k$ pointwise.
\end{lem}
\begin{proof}
It is true for $k = 0.$  To compute $v_{k+1}$ one solves the
differential equation by the method of variation of
parameters,\footnote{We can write $G_0 h(t) = \int_T^\infty
K_0(t,s)h(s)ds,$ where
$K_0(t,s)=\frac{1}{3}(\frac{s^2}{t}-\frac{t^2}{s})$ for $s>t$ and
zero otherwise (this form of the kernel is relevant to the
condition of good functional behavior at infinity).  One can also
use the equivalent $G_0 h(t) = t^2 \int_t^\infty
t_1^{-4}\int_{t_1}^\infty t_2^2 h(t_2)dt_2dt_1,$ which appears in
the text.} using the homogeneous solution $t^2$ or $t^{-1}.$ We
have
$$v_{k+1}(t)= t^2\int_{t}^{\infty}t_1^{-4}\int_{t_1}^\infty t_2^2
(-4t_2e^{-t_2}) v_k(t_2)dt_2dt_1.$$ One computes $v_1(t)=-4(t + 2 +
2/t)e^{-t},$ and therefore $|v_1(t)|<16te^{-t}$ for $t>2$.  Now
assume $|v_k(t)|<(16te^{-t})^k$ for some $k\geq 1.$  First compute
for $a,b\in \mathbb{N},$ $t> a/b$ and $t>2$:
\begin{eqnarray*}
\int_{t}^\infty s^a e^{-bs}ds&=&-\frac{1}{b}\left.\left[s^a +
\frac{a}{b}s^{a-1}+
\frac{a(a-1)}{b^2}s^{a-2}+...\right]e^{-bs}\right|_{t}^\infty\\
&\leq&\frac{1}{b}t^a(1 + (a/bt) + (a/bt)^2 + ...)e^{-bt}\\
&\leq&\frac{t}{b(t-1)}t^a e^{-bt}\leq 2t^ae^{-bt}.
\end{eqnarray*}
Therefore,
\begin{eqnarray*}
|v_{k+1}(t)|&\leq& 4t^2\int_t^\infty t_1^{-4}\int_{t_1}^\infty t_2^3 e^{-t_2}v_k(t_2)dt_2dt_1\\
&\leq& 4t^2\int_t^\infty t_1^{-4}\int_{t_1}^\infty t_2^3 e^{-t_2}(16)^k t_2^k e^{-kt_2}dt_2dt_1\\
&\leq& (16)^k\cdot 4\cdot 2t^2\int_t^\infty t_1^{k+3-4}e^{-(k+1)t_1}dt_1\\
&\leq& (16)^{k+1}t^{k+1}e^{-(k+1)t}
\end{eqnarray*}
and the lemma is proven.

Note that we needed $k-1\geq 0$ to bound $a/b$ from above. That the
$k=1$ term is the proper order follows from some fortuitous
cancelation.
\end{proof}

It now follows that $u < C(1-4te^{-t})^{-1} 4te^{-t} \leq C'
te^{-t}$ for some constant $C'.$  The space $S=\{g(t):|g(t)|\leq 2C'
te^{-t}\}$ forms a closed subset of $\mathcal B$ on which we apply
the contraction mapping theorem.

The proof of the previous lemma also shows the following.

\begin{lem} \label{quad-bound}
If $|h(t)| \le C(te^{-t})^2$, then there are positive constants $T$
and $K$ independent of $h$ so that if $t\ge T$, then
$$|(Gh)(t)| \le CK(te^{-t})^2.$$
\end{lem}

\begin{proof}
The computations above show that if $\ell\ge 2$ and $|w(t)|\le C''
(te^{-t})^\ell$, then $|(L_0^{-1}w)(t)| \le 4C''(te^{-t})^\ell$.
Then for our $h(t)$, compute
\begin{eqnarray*}
|(Gh)(t)| &\le& |(L_0^{-1}h)(t)| + |(L^{-1}_0 \delta_L L_0^{-1}
h)(t)| + \cdots \\
&\le& 4C(te^{-t})^2 + 4\cdot16(te^{-t})^3 + 4\cdot 16^2 (te^{-t})^4
+ \cdots \\
&=& \frac{4C(te^{-t})^2}{1-16te^{-t}},
\end{eqnarray*}
and so for $T$ large enough, we can choose $K=4/(1-16Te^{-T})$.
\end{proof}

\begin{prop} There is a constant $T>0$ so that for $t\ge T$,
the equation (\ref{psi-eq}) has a solution of the form $\log(2t) + t
+ O(te^{-t}).$
\end{prop}
\begin{proof}  We now show the mapping $\phi\rightarrow A\phi \equiv
Gf - GQ\phi$ is a contraction. First, $Gf$ lies within $S$ and
$Q\phi$ is small since $Q$ is a quadratic nondifferential operator.
More specifically, for some $T>0$, the sup norm of $\phi$ on
$[T,\infty)$ can be made arbitrarily small.  Therefore, by Lemma
\ref{quad-bound}, $\|GQ\phi\|_{\mathcal B} \ll \|\phi\|_{\mathcal
B}$ on $[T,\infty)$. As a result, $A$ maps $S$ to $S.$ Further, note
$\|A\phi_1 - A\phi_2\|_{\mathcal B} = \|GQ\phi_1 -
GQ\phi_2\|_{\mathcal B} = \|G(Q\phi_1-Q\phi_2)\|_{\mathcal B},$
since $G$ is linear.

Since $Q$ is quadratic and $\phi_1,\phi_2\in S$, $$|(Q\phi_1 -
Q\phi_2)(t)| \le K'(te^{-t})^2\|\phi_1-\phi_2\|_{\mathcal B},$$ for
$K'$ depending on $C'$. Then Lemma \ref{quad-bound} shows
$$|(GQ\phi_1 - GQ\phi_2)(t)| \le KK'(te^{-t})^2
\|\phi_1-\phi_2\|_{\mathcal B}.$$ Then clearly, by taking $T$ large
enough, there is a fixed $\theta<1$ so that
$\|A\phi_1-A\phi_2\|_{\mathcal B} \le \theta \|\phi_1 -
\phi_2\|_{\mathcal B}$ for $\phi_1,\phi_2\in S$. By the fixed point
theorem, there exists $\phi$ such that $A\phi = \phi.$ $\phi$ is
smooth by standard bootstrapping. Then $\log(2t) + t + \phi(t)$
solves (\ref{psi-eq}).
\end{proof}

In the next section, we will make an ansatz that the local form of
our global function be consistent with the dominant $\log|\log|z|^2|
- \log|z|$ behavior of this local solution. The dominant term is
$-\log|z|$ which comes from the form of $U$ and determines the
residue at the singularity.

\subsection{Global Existence}

The coordinate-independent version of equation (\ref{psi-eq}) for a
general background metric is
\begin{equation}
\label{eqone} \Delta u + 4\Vert U\Vert^2 e^{-2u} + 2e^u - 2\kappa
= 0
\end{equation}
on $S^2,$ where norms, the Laplacian, integrals, and the Gauss
curvature $\kappa$ are taken with respect to the background metric.
$U$ is a holomorphic cubic differential, which we take to have
exactly 3 poles of order 2 and thus no zeroes, and $u$ is taken to
have a prescribed singularity structure such that $\int \Delta u =
6\pi,$ which follows from our local analysis in Section \ref{local}.

Near each pole of $U$, there is a local coordinate $z$ so that the
pole is at $z=0$, and $U=z^{-2}dz^3$ exactly.  We call this $z$
the \emph{canonical holomorphic coordinate}.  In a neighborhood of
each pole, we take
\begin{equation}
u_0 = \log|\log|z|^2| - \log|z|
 \end{equation}
and the background metric to be $|dz|^2$.  This $u_0$ is an explicit
solution of the parabolic affine sphere equation (\ref{inteq}). The
background metric and $u_0$ are extended smoothly to the rest of
$\cp^1$.  Note that $\int \Delta u_0=6\pi$ (each pole contributes
$2\pi$).  All integrals in this section will be evaluated with
respect to the background metric.

To implement the required singularity structure, we write $u = u_0
+ \eta$ for $\eta$ in the Sobolev space $H_1$.  Note this implies
$\int\Delta\eta=0$.  We define the functional
\begin{equation}
\label{Jeq}
\begin{array}{rcl}
J(\eta) &=& \int \left( \frac{1}{2}|\nabla \eta|^2 + (2\kappa
-\Delta u_0)\eta + \frac{1}{2}3\cdot 4\Vert U\Vert^2
e^{-2u_0}e^{-2\eta}\right) \\ &\phantom{=}& \\
 &\phantom{=}& - 2\pi \log \int
\left( 4\Vert U\Vert^2 e^{-2u_0}e^{-2\eta} + 2 e^{u_0}e^\eta
\right).
\end{array}
\end{equation}
(We note that it is necessary to separate $\eta$ from $u_0$, as
$\nabla u_0$ is not in $L^2.$)  $J$ is not defined for all functions
$\eta\in H_1$.  One problem is that $\Delta u_0\notin L^2$. The term
$\int\Delta u_0\eta$ can be taken care of by integrating by parts
(see the proof of Proposition \ref{J-bound-prop} below). A more
serious problem is that $4\|U\|^2e^{-2u_0}\notin L^p$ for any $p>1$.
This cannot be fixed by integrating by parts, as the example $\eta =
-\frac12 \log|\log|z|^2|\in H_{1,{\rm loc}}$ shows. That said, there
is a uniform lower bound on $J$ among all $\eta\in H_1$ so that
$\int 4\|U\|^2e^{-2u_0} e^{-2\eta}<\infty$ (see the remark after
Proposition \ref{J-bound-prop}). Thus we can still talk of taking
sequences of $\eta\in H_1$ to minimize $J$. (The term
$\int2e^{u_0}e^{\eta}$ is always finite for $\eta\in H_1$ since
$e^{u_0}\in L^p$ for $p<2$ and Moser-Trudinger shows that $e^\eta\in
L^q$ for all $q<\infty$.)

We wish to show that $J(\eta)$ has a local minimum. If so, then the
minimizer satisfies the Euler-Lagrange equation
 \begin{equation}\label{el-frac}
 \Delta \eta - (2\kappa -
\Delta u_0) + 3\cdot 4\Vert U\Vert^2 e^{-2u_0}e^{-2\eta} +
\frac{-2\cdot 4\Vert U\Vert^2 e^{-2u_0}e^{-2\eta} + 2
  e^{u_0}e^\eta}{\frac{1}{2\pi}\int 4\Vert U\Vert^2e^{-2u_0}e^{-2\eta} + 2
  e^{u_0}e^\eta} = 0.
  \end{equation}
One can easily check by integrating this equation that for a
solution $\eta_0,$ the denominator in the last term must be equal
to one.  Thus $u = \eta_0 + u_0$ satisfies the original equation
(\ref{eqone}).  In this case, the equation $\eta_0$ satisfies is
\begin{equation} \label{eta-sol}
\Delta \eta - (2\kappa - \Delta u_0) + 4\Vert U\Vert^2
e^{-2u_0}e^{-2\eta} + 2 e^{u_0}e^\eta = 0.
\end{equation}
This is equivalent to equation (\ref{psi-eq}), the equation for the
metric of an elliptic affine sphere: for the background metric $h$,
write $e^{u_0+\eta}h=e^\psi|dz|^2$. Then $\eta$ satisfies
(\ref{eta-sol}) if and only if $\psi$ satisfies (\ref{psi-eq}).

\begin{dfn}
We call $\eta$ \emph{admissible} if $\eta \in H_1$ and $\int
4\Vert U\Vert^2 e^{-2u_0}e^{-2\eta} <\infty.$
\end{dfn}

In order to analyze the functional $J$, for an admissible $\eta$,
consider $J(\eta+k)$ for $k$ a constant.  $J(\eta+k)$ has the form
$$(\text{indep.\ of }k) + 2\pi k +
3\pi Ae^{-2k}-2\pi\log[2\pi(Ae^{-2k} +Be^k)],$$
 where
 \begin{equation}\label{defA-B}
 A = A(\eta) \equiv \frac{1}{2\pi}\int
4\Vert U\Vert^2 e^{-2u_0}e^{-2\eta}, \quad B = B(\eta) \equiv
\frac{1}{2\pi}\int 2e^{u_0}e^\eta.
 \end{equation}
Thus upon setting $(\partial/\partial k)J(\eta+k)=0$, we find a
critical point only if $Ae^{-2k} + Be^k = 1$, and this can only
happen if
$$AB^2\le \frac4{27}.$$
 If $AB^2>4/27$, then the infimum occurs as $k\to+\infty$, and if
 $AB^2<4/27$, there are two finite critical points: a local
 minimum for which $B(\eta+k)=Be^k<2/3$ and a local maximum for
 which $B(\eta+k)>2/3$. With that in mind, we formulate the
 following variational problem:

Let $$Q = \{\eta \in H_1 : A + B \le 1\}.$$ We will minimize $J$
for $\eta \in Q$.  Note that this will avoid the potential problem
at $k\to+\infty$, where $B(\eta+k)\to+\infty$.  Also, the
inequality in the definition of $Q$ will be important. It will
allow us to use the Kuhn-Tucker conditions to control the sign of
the Lagrange multiplier in the Euler-Lagrange equations.  The
discussion above about adding a constant $k$ can be summarized in

\begin{lem} \label{min-on-line}
If $\eta\in Q$, then the minimizer of $$\{J(\eta+k): k \mbox{
constant, }\eta+k\in Q\}$$
 occurs for $k$ so that $A(\eta+k)+B(\eta+k)=1$, $B(\eta+k)\le
2/3$, and $k\le0$. If $A(\eta)+B(\eta)<1$, then the minimizer
$k<0$ and $B(\eta+k)<2/3$. Moreover, if $A(\eta)+B(\eta)=1$ and
$B(\eta)\le2/3$, then $k=0$.
\end{lem}
\begin{proof} Compute $(\partial/\partial k)J(\eta+k)$ and use the
first derivative test.
\end{proof}

\begin{prop} \label{J-bound-prop}
There are positive constants $\gamma$ and $R$ so that for all
$\eta \in Q$,
$$ J(\eta) \ge \gamma \int |\na\eta|^2 - R.$$
\end{prop}

\begin{rem}
We can also prove the same result for all admissible $\eta\in
H_1$.  In this case, we must also control potential minimizers at
$k=+\infty$.  For admissible $\rho\in H_1$ so that $\int\rho=0$,
consider the functional $$\tilde J(\rho)=\lim_{k\to\infty}
J(\rho+k).$$ We bound $\tilde J$ from below much the same as the
following argument, although there also is an extra term in
$\tilde J$ that must be handled using the Moser-Trudinger
estimate.
\end{rem}

\begin{proof}
As above, $u_0=\log|\log|z|^2| - \log|z|$ in the canonical
coordinate $z$ near each pole of $U$.  Since $\Delta u_0 \notin
L^2$, we should integrate by parts to handle the $-\int \Delta u_0
\,\eta$ term in $J$.  Let $u_0'=\log|\log|z|^2|$ near each pole of
$U$ and smooth elsewhere.  Then $\Delta u_0 = \Delta u_0'$ near
each pole and the difference $\Delta u_0 -\Delta u_0'$ is smooth
on $\cp^1$.  Then if we let $\zeta$ be the smooth function
$2\kappa -\Delta (u_0-u_0')$,
\begin{eqnarray*}
J(\eta) &=& \int [\frac12|\na\eta|^2 + \zeta \eta -\Delta
u_0'\,\eta] + 3\pi A- 2\pi \log 2\pi(A+B)\\
&>& \int [\frac12|\na\eta|^2 + \zeta \eta +\na u_0'\cdot \na \eta]
- 2\pi \log 2\pi \\
&\ge& C + \int \left[\frac12|\na\eta|^2 - \frac1{4\epsilon}\zeta^2
-\epsilon\eta^2 - \frac1{4\epsilon}|\na u_0'|^2 - \epsilon|\na
\eta|^2\right] \\
&\ge& C_\epsilon + \int \left(\frac12-\delta\right) |\na\eta|^2
\\
\end{eqnarray*}
 Here $\delta = (\frac1{\lambda_1}+1)\epsilon$, for $\lambda_1$
the first nonzero eigenvalue of the Laplacian of the background
metric, and we've used the facts that $A>0$ and $A+B\le1$.
\end{proof}

Here is another useful lemma.
\begin{lem} \label{ABsquared} For any $\eta\in H_1$,
$$ A B^2 \ge L=2\pi^{-3} \left(\int
\|U\|^\frac23 \right)^3.$$
 If $AB^2 = L$, then there is a constant $C$ such that
 $$\eta = C +\sfrac23 \log\|U\| - u_0.$$
\end{lem}

\begin{proof}
Let $f=(4\|U\|^2)^\frac13 e^{-\frac23(u_0+\eta)}$,
$g=e^{\frac23(u_0+\eta)}$. Apply H\"older's inequality $\int fg
\le \|f\|_3 \|g\|_{\frac32}$. The last statement follows from the
case of equality in H\"older's inequality.
\end{proof}

\begin{rem}
The bound $L$ in the previous lemma does not depend on the
background metric; it depends only on the conformal structure on
$\cp^1$ and the cubic form $U$.
\end{rem}

An admissible $\eta\in H_1$ is a \emph{weak solution} of
(\ref{eta-sol}) if $\eta$ is a solution of (\ref{eta-sol}) in the
sense of distributions.

\begin{prop} \label{min-solve}
Assume that $U$ is such that $L<4/27$. Then any minimizer $\eta$
of $\{J(\eta):\eta \in Q\}$ is a weak solution of (\ref{eta-sol}).
\end{prop}
\begin{proof}  Recall $Q=\{\eta:A+B\le1\}$.

Case 1: The minimizer $\eta$ satisfies $A+B<1$. Since the
constraint $A+B\le 1$ is slack, $\eta$ must satisfy the
Euler-Lagrange equation (\ref{el-frac}). Then as above, we may
integrate to find that the denominator $A+B$ in (\ref{el-frac})
must be equal to $1$. Thus this case cannot occur.

Case 2: The minimizer $\eta$ satisfies $A+B=1$. In this case, we
have Lagrange multipliers $[\mu_0,\mu_1]\in\rp^1$ so that $\eta$
weakly satisfies
$$\begin{array}{c}
\D \mu_0\left[\Delta \eta - (2\kappa-\Delta u_0) + 3\cdot4\|U\|^2
e^{-2u_0} e^{-2\eta} + \frac{-2\cdot4\|U\|^2e^{-2u_0}e^{-2\eta} +
2e^{u_0}e^\eta}{A+B}
\right] \\ \ \\
= \mu_1(-2\cdot4\|U\|^2e^{-2u_0}e^{-2\eta} + 2e^{u_0}e^\eta),
\end{array}$$
and $A+B=1$.  Thus,
 \begin{equation}\label{lag-mult}
 \mu_0[\Delta \eta - (2\kappa-\Delta u_0)
+ a + b] = \mu_1(-2a+b)
\end{equation}
 for $$a =4\|U\|^2e^{-2u_0}e^{-2\eta}, \qquad
 b=2e^{u_0}e^\eta.$$
 Note then that $A=\int a/2\pi$, $B=\int b/2\pi$.

Also note the constraint the Kuhn-Tucker conditions place on the
Lagrange multipliers.  Recall that if we minimize a function $f$
subject to the constraint $g\le1$, and if the minimum occurs on
the boundary $g=1$, then we have $\mu_0\na f=\mu_1 \na g$ for
$\mu_0\mu_1\le 0$. This is exactly our situation for $f=J$ and
$g=A+B$.

Thus we have three cases: if $\mu_1=0$, then equation
(\ref{lag-mult}) becomes equation (\ref{eta-sol}) and we've proved
the proposition.

In the second case, if $\mu_0=0$, then the Euler-Lagrange equation
(\ref{lag-mult}) may be solved explicitly for $\eta$ to find
$$ \eta = \sfrac13 \log (4\|U\|^2) - u_0.$$
Near each pole of $U$, there is a coordinate $z$ so that $\|U\| =
|z|^{-2}$ and $u_0 = \log|\log|z|^2|-\log|z|$.  So $$\eta =
\sfrac13\log4 -\sfrac13\log|z| - \log|\log|z|^2|$$ there and so
$\eta\notin H_1$.

Finally, we consider where $\mu = \mu_1/\mu_0<0$.  We will analyze
the second variation at any critical point to show that there are
no minimizers in this case.

Integrate (\ref{lag-mult}) to find
 $$-2\pi + 2\pi A + 2\pi B = \mu(-2\cdot2\pi A +
 2\pi B).$$
 Then since $A+B=1$, we have $2A=B$, since we are in the case $\mu\neq0$.
So $A=1/3$ and $B=2/3$.  We analyze the second variation to show
that for $L<4/27$, there is no minimizer at $A=1/3$, $B=2/3$
(unless possibly if $\mu_1=0$).

Let $\eta$ satisfy (\ref{lag-mult}) and $A=1/3,$ $B=2/3$. Consider
a variation $\eta + \epsilon \alpha + \frac{\epsilon^2}2\beta$ so
that $\eta$ satisfies $A+B=1$ to second order when
$\epsilon=0$.\footnote{This corresponds to an actual variation in
$Q$ by standard Implicit Function Theorem arguments---see
\cite{lang85}. Let $X$ be the Banach space $H_1\cap C^0$.  Then
let $g\!:X\to\re$, $g(\nu)=A(\eta+\nu)+B(\eta+\nu)$.  It is
straightforward to show that $g$ is $C^1$ in the Banach space
sense.  Moreover, for $2a\neq b$ (which holds for any $\eta\in
H_1$), we can check that $dg\!:X\to\re$ is nonzero.  So then
$Y=g^{-1}(1)=\{A+B=1\}$ is a Banach submanifold of $X$ near
$\nu=0$. So for any element $\alpha\in\ker dg_0$, there is a curve
in $Y$ tangent to $\alpha$. Along such a curve, we compute
restrictions on the second-order term $\beta$.}

We assume $\alpha$ is a constant.  Then the first variation
$$\left.\frac{\partial}{\partial\epsilon}(A+B)\right|_{\epsilon=0}
= -2\alpha A + \alpha B = 0$$ for $A=1/3$, $B=2/3$.  So to first
order $\eta+\alpha$ satisfies $A+B=1$ and $\alpha$ is tangent to
$\{A+B=1\}$.

Now we require
\begin{eqnarray}
 0&=& 2\pi\left.\frac{\partial^2}{\partial\epsilon^2}(A+B)\right|_
{\epsilon=0} \nonumber \\
&=& \int a(4\alpha^2-2\beta) + b(\alpha^2+\beta) \nonumber \\
&=& \alpha^22\pi(4A+B) + \int\beta(-2a+b)  \nonumber \\
&=& 2\pi\cdot2\alpha^2 + \int\beta(-2a+b) \label{beta-constr}.
\end{eqnarray}

Now for this variation
$J=J(\eta+\epsilon\alpha+\frac{\epsilon^2}2\beta)$, compute
\begin{eqnarray}
\left.\frac{\partial^2J}{\partial\epsilon^2}\right|_{\epsilon=0}
&=& \int \na\eta\cdot\na\beta + |\na\alpha|^2 + (2\kappa-\Delta
u_0)\beta + \sfrac32 \int a(4\alpha^2-2\beta)  \nonumber \\
&& {}- 2\pi \frac{\int a(4\alpha^2-2\beta)+b(\alpha^2+\beta)}
{\int a + b} + 2\pi \frac{\left(\int a(-2\alpha)+b\alpha\right)^2}
{\left(\int a+b\right)^2} \nonumber \\
&=& \int [\na\eta\cdot\na\beta + (2\kappa-\Delta u_0)\beta -
3a\beta] + 2\pi\cdot 6\alpha^2 A \nonumber \\
&&{}-2\pi\cdot \alpha^2(4A+B) - \int\beta(-2a+b) \nonumber \\
&=& \int [\na\eta\cdot\na\beta + (2\kappa-\Delta u_0)\beta -
3a\beta] + 2\pi\cdot 2\alpha^2  \label{second-var}
\end{eqnarray}
Here we've used the following facts to get from the first line to
the second: $\na\alpha=0$ since $\alpha$ is constant, $\int
(a+b)/2\pi = A+B=1$, and the last term vanishes since $\alpha$ is
constant and $2A=B$.  The third line follows from the second by
the constraint (\ref{beta-constr}) and the fact $A=1/3$.

Now we use the Euler-Lagrange equation (\ref{lag-mult}).  Recall
$\mu_0\neq0$ and $\mu=\mu_1/\mu_0$. Then
$$\int\na\eta\cdot\na\beta = -\int(\Delta\eta)\beta = \int[-(2\kappa-
\Delta u_0) + a + b - \mu(-2a+b)]\beta.$$
 Plug this into (\ref{second-var}) to find
 \begin{eqnarray*}
\left.\frac{\partial^2J}{\partial\epsilon^2}\right|_{\epsilon=0}
&=& (1-\mu)\int(-2a+b)\beta +2\pi\cdot 2\alpha^2 \\
 &=& 2\pi\cdot2\mu\alpha^2.
\end{eqnarray*}
Here the last line follows from (\ref{beta-constr}).  Thus if we
choose $\alpha\neq0$, then the second variation along this path is
negative since $\mu<0$.  Therefore, there is no minimizer for our
variational problem satisfying $\mu<0$.
\end{proof}

Now we show that there is a minimizer.

\begin{lem} \label{AB-bound}
Assume $L<4/27$. Then there is a constant $\delta>0$ so that $A,B
\in(\delta,1/\delta)$ for all $\eta\in Q$.
\end{lem}
\begin{proof}
Lemma \ref{ABsquared} implies that $AB^2\ge L$. Since $0<L<4/27$,
$A>0$, $B>0$, and $A+B\le 1$, this proves the lemma.
\end{proof}

\begin{lem} \label{jensen}
There are constants $K_1$, $K_2$ so that for all admissible
$\eta\in H_1$, and for $c=(\int\eta)/(\int1)$,
$$\log A \ge K_1 -2c, \qquad
\log B \ge K_2+c.$$
\end{lem}
\begin{proof}
Since exp is convex, Jensen's inequality gives
\begin{eqnarray*}
\log A &=& -\log 2\pi + \log \int 4\|U\|^2 e^{-2u_0}
e^{-2\eta} \\
&\ge& -\log2\pi  + \frac{\int\log4\|U\|^2 -2u_0 - 2\eta}{\int 1} +
\log \int1.
\end{eqnarray*}
The case for $B$ is the same.
\end{proof}

\begin{lem}
Let $\eta_i$ be a sequence in $Q$ so that $\lim_i J(\eta_i) =
\inf_{\eta\in Q} J(\eta).$ Then there is a positive constant $C$
so that $\|\eta_i\|_{H_1} \le C$ for all $i$.
\end{lem}
\begin{proof}
First we note that Lemmas \ref{AB-bound} and \ref{jensen} show
that the average value $c=(\int\eta)/(\int 1)$ is uniformly
bounded above and below for all $\eta\in Q$.

Proposition \ref{J-bound-prop} shows that $J(\eta)\ge \gamma
\int|\na\eta|^2 -R$ for $\gamma, R>0$ uniform constants. Thus for
any minimizing sequence, $\int|\na\eta|^2$ must be uniformly
bounded. Then write $\eta=\rho+c$ for $\int\rho=0$, $c$ constant.
Then
$$\|\eta\|_{L^2} \le \|\rho\|_{L^2} + \|c\|_{L^2} \le
\lambda_1^{-\frac12}\|\na\rho\|_{L_2} + K =
\lambda_1^{-\frac12}\|\na\eta\|_{L_2} + K$$
 for $K$ a uniform constant and $\lambda_1$ the first nonzero
 eigenvalue of the Laplacian.  This shows the $H_1$ norm of $\eta$
 in the minimizing sequence is uniformly bounded.
\end{proof}

Now given a minimizing sequence $\{\eta_i\} \subset Q$, Lemma
\ref{min-on-line} shows that we can assume
$A(\eta_i)+B(\eta_i)=1$, $B(\eta_i)\le2/3$. Then there is a
subsequence, which we still refer to as $\eta_i$,  which is weakly
convergent to a function $\eta_\infty \in H_1$ (the weak
compactness of the unit ball in a Hilbert space), strongly
convergent to $\eta_\infty$ in $L^p$ for $p<\infty$ (Sobolev
embedding), convergent pointwise almost everywhere to
$\eta_\infty$ ($L^p$ convergence implies subsequential
almost-everywhere convergence), and so that $e^{\eta_i}$ is
strongly convergent to $e^{\eta_\infty}$ in $L^p$ for $p<\infty$
(Moser-Trudinger). Recall
$$J(\eta)=\int[\sfrac12|\na\eta|^2 + (2\kappa -\Delta u_0)\eta]
+ 3\pi A - 2\pi \log 2\pi(A+B).$$
 Then the second term in the integral converges by strong
 convergence in $L^1$ and weak convergence in $H_1$ (see the proof
of Proposition \ref{J-bound-prop} for the integration by parts
trick). The term $\int\sfrac12|\na\eta|^2$ is lower semicontinuous
(the norm in a
 Hilbert space is lower semicontinuous under weak convergence).
 Lower semicontinuity is enough since we are seeking a minimizer.
 $B$ converges by Moser-Trudinger: $e^{u_0}\in L^p$ for $p<2$.
 Then since $e^{\eta_i}$ converges in $L^q$ for
 $\frac1p+\frac1q=1$, $B=\int2e^{u_0}e^\eta$ converges.

That leaves the term $A$.  Fatou's lemma and the almost-everywhere
convergence of $\eta_i$ then show
$$ A(\eta_\infty) \le \liminf_{i\to\infty} A(\eta_i).$$
We want to rule out the case of strict inequality. Note
$A(\eta_\infty) + B(\eta_\infty) \le \lim A(\eta_i)+B(\eta_i) =
1$, and so $\eta_\infty\in Q$. Also, since $A(\eta_i)+B(\eta_i)=1$
and $B(\eta_i)\to B(\eta_\infty)$, $\lim A(\eta_i) =
1-B(\eta_\infty)$.

Consider the constant $k$ so that $\eta_\infty + k$ minimizes
$$ \{J(\eta_\infty + k) : \eta_\infty+k \in Q\}.$$
Note that Lemma \ref{min-on-line} shows that
$e^{-2k}A(\eta_\infty) + e^k B(\eta_\infty)=1$. Now compute

\begin{eqnarray*}
\lim_{i\to\infty} J(\eta_i) &\ge&
\int\left[\,\sfrac12|\na\eta_\infty|^2+(2\kappa-\Delta
u_0)\eta_\infty\right] +
3\pi\left[1-B(\eta_\infty)\right], \\
J(\eta_\infty+k)&=&
\int\left[\,\sfrac12|\na\eta_\infty|^2+(2\kappa-\Delta
u_0)(\eta_\infty+k)\right] + 3\pi e^{-2k}A(\eta_\infty). \\
\end{eqnarray*}
Now substitute $e^{-2k}A(\eta_\infty) = 1-e^kB(\eta_\infty)$ to
show
\begin{equation} \label{lim-dif}
\lim_{i\to\infty} J(\eta_i)-J(\eta_\infty+k) \ge -2\pi k
+3\pi B(\eta_\infty)(e^k-1)\\
\end{equation}

We prove $A(\eta_\infty)=\lim A(\eta_i)$ by contradiction.  If on
the contrary $A(\eta_\infty)<\lim A(\eta_i)$, Lemma
\ref{min-on-line} and the fact $A(\eta_i)+B(\eta_i)=1$ imply that
$k<0$. Then it is straightforward to check that the right-hand side
of (\ref{lim-dif}) is strictly positive (it is zero if $k=0$, and
its derivative with respect to $k$ is negative for $k<0$---use the
fact $B(\eta_\infty)=\lim B(\eta_i)\le 2/3$). This shows $\lim
J(\eta_i) > J(\eta_\infty+k)$ and so contradicts the fact that
$\eta_i$ is a minimizing sequence for $J$.

The same analysis shows that $\lim \int|\na\eta_i|^2 =
\int|\na\eta_\infty|^2$. So $J(\eta_\infty)=\lim J(\eta_i)$, and
$\eta_\infty$ is a minimizer of $\{J(\eta):\eta\in Q\}$.

\begin{thm} If $L<4/27$ then a weak solution to (\ref{eta-sol})
exists.  Conversely, if $L\ge4/27$, then there is no weak solution
to (\ref{eta-sol}).
\end{thm}
\begin{proof}
The preceding paragraphs, together with Proposition
\ref{min-solve}, prove existence in the case $L<4/27$. We address
the nonexistence in two cases:

Case  $L>4/27$. If $\eta$ solves (\ref{eta-sol}), then we can
integrate (\ref{eta-sol}) to find $A+B=1$. On the other hand,
$A>0$, $B>0$, and Lemma \ref{ABsquared} shows that $AB^2\ge L >
4/27$. Simple calculus shows that there is no such pair $(A,B)$ in
this case.

Case $L=4/27$. As in Case 1, we must have $A+B=1$ and $AB^2\ge
L=4/27$. The only way this can happen is if $A=1/3$, $B=2/3$, so
that $AB^2=4/27$.  In this case, Lemma \ref{ABsquared} forces
$\eta = C +\frac23 \log\|U\|-u_0$ for some constant $C$. Since
$u_0 = \log|\log|z|^2| -\log|z|$ and $\|U\|=|z|^{-2}$ near each
pole of $U$, $\eta = C -\log|\log|z|^2| - \frac13\log|z|$ near
each pole of $U$. Thus $\eta\notin H_1$.
\end{proof}

\begin{prop} Any weak solution $\eta$ to (\ref{eta-sol}) is smooth away
from the poles of $U$.
\end{prop}
\begin{proof}
In a neighborhood bounded away from the poles of $U$, the
quantities $\|U\|^2$ and $u_0$ are smooth and bounded.  Since
$\eta\in H_1$, Moser-Trudinger shows that $e^\eta,\,e^{-2\eta} \in
L^p$ for all $p<\infty$.  Therefore, (\ref{eta-sol}) implies
$\Delta \eta \in L^p_{\rm loc}$.  Since $\eta\in L^p$ by Sobolev
embedding, the $L^p$ elliptic theory \cite{gilbarg-trudinger}
shows that $\eta\in W^{2,p}_{\rm loc}$. Sobolev embedding shows
$\eta\in C^{0,\alpha}_{\rm loc}$, and so $\Delta\eta\in
C^{0,\alpha}_{\rm loc}$. The Schauder theory then shows $\eta\in
C^{2,\alpha}_{\rm loc}$. Further bootstrapping implies $\eta$ is
smooth.
\end{proof}

\begin{rem}
It is not clear whether the solution constructed is unique.  The
maximum principle does not work to give uniqueness.
\end{rem}

\subsection{A metric for the ``Y'' vertex}

Let $\widetilde{\Sigma}$ be the universal cover of ${\Sigma} =
S^2\setminus \{p_1,p_2,p_3\}.$  Lifting the appropriate objects to
the cover we find a solution to (\ref{eqone}) on
$\widetilde{\Sigma}.$ Since the equation (\ref{eqone}) is the
integrability condition for the developing map, we have a solution
$\widetilde{f}:\widetilde{\Sigma}\rightarrow \R^3,$ with monodromies
of $\Sigma$ acting as equiaffine deck transformations fixing the
normal vector $\xi$ and acting by isometry. The quotient by the deck
transformations gives an elliptic affine sphere structure on
$\Sigma$ as well as the locally defined developing map $f.$ Then the
map $F:(\Sigma \times \R_+)\rightarrow \R^3$ defined by $F(x,r)=
rf(x) =: (y_1, y_2, y_3)$ maps the cone over $\Sigma$ to $\R^3$ and
is locally invertible (so we may express $r = r(y)$). The potential
function $\Phi(y) = r^2/2$ defines a parabolic affine sphere on a
neighborhood of the ``Y'' vertex, by Baues and Cort\'es's Theorem
\ref{bc-thm}.  This is our main result.

\begin{rem}
The monodromy group of this metric determines the affine flat
structure.  We have not yet determined this monodromy group, thus
cannot verify that the metric is one predicted by Gross-Siebert
\cite{gs} and Haase-Zharkov \cite{hz03}.
\end{rem}

\vskip 0.2in \centerline{\bf Acknowledgments} \vskip 0.1in We would
like to thank Rafe Mazzeo for stimulating discussions, and the
referees for very useful comments and clarifications, in particular
the link to the work of McIntosh. The work of S.-T. Yau was
supported in part by the National Science Foundation (DMS-0244464,
DMS-0074328, DMS-0306600, DMS-9803347). The work of E. Z. was
supported in part by the National Science Foundation (DMS-0072504)
and by the Alfred P. Sloan Foundation.

\bibliographystyle{abbrv}

\begin{thebibliography}{99}

\bibitem{AKMV} M. Aganagic, A. Klemm, M. Mari\~no, and C. Vafa,
``The Topological Vertex,'' Comm. Math. Phys. {\bf 254} (2005)
425--478, hep-th/0305132, MR2117633.

\bibitem{amari-nagaoka} S. Amari and H. Nagaoka, {\sl Methods of
Information Geometry,} The American Mathematical Society,
Providence, 2000, MR1800071, Zbl 0960.62005.

\bibitem{BC} O. Baues and V. Cort\'es,
``Proper Affine Hyperspheres which Fiber over Projective Special
Kaehler Manifolds,'' Asian J. Math. {\bf 7} (2003) 115--132,
math.DG/0205308, MR2015244.

\bibitem{calabi3} E. Calabi, ``Affine Differential Geometry and
Holomorphic Curves," in {\sl Complex Geometry and Analysis (Pisa,
1988)}, Lecture Notes in Mathematics {\bf 1422}, Springer-Verlag,
(1990) 15--21, MR1055839, Zbl 0697.53019.

\bibitem{calabi2} E. Calabi, ``Complete Affine Hypersurfaces I,''
in {\sl Symposia Mathematica {\bf X}}, Academic Press, London
(1972) 19--38, MR0365607, Zbl 0252.53008.

\bibitem{calabi} E. Calabi, ``A Construction of
Nonhomogeneous Einstein Metrics,'' Proc. of Symp. in Pure
Mathematics {\bf 27}, AMS, Providence (1975) 17--24, MR0379912,
Zbl 0309.53043.

\bibitem{cheng-yau86}
S.-Y. Cheng and S.-T. Yau.
\newblock Complete affine hyperspheres. part {I}. {T}he completeness of affine
  metrics.
\newblock {\em Communications on Pure and Applied Mathematics}, 39(6):839--866,
  1986, MR0859275, Zbl 0623.53002.

\bibitem{cortes} V. Cort\'es, ``A Holomorphic Representation Formula
for Parabolic Hyperspheres,'' in {\sl PDEs, Submanifolds, and
Affine Differential Geometry (Warsaw, 2000)} Banach Center Publ.
{\bf 57}, Polish Acad. Sci., Warsaw (2002) 11--16,
math.DG/0107037, MR1972459, Zbl 1029.53017.

\bibitem{ferrermm} L. Ferrer, A. Mart\'inez, and F. Mil\'an, ``An
Extension of a Theorem by K. J\"orgens and a Maximum Principle at
Infinity for Parabolic Affine Spheres," Math. Z. {\bf 230} (1999)
471--486, MR1679973, Zbl 0967.53009.

\bibitem{Fuk} K. Fukaya, ``Multivalued Morse Theory,
Asymptotic Analysis, and Mirror Symmetry,'' in {\sl Graphs and
Patterns in Mathematics and Theoretical Physics} Proc. Sympos.
Pure Math. {\bf 73}, AMS, Providence, RI (2005) 205--278,
MR2131017.

\bibitem{gilbarg-trudinger} D. Gilbarg and N. Trudinger, {\sl Elliptic Partial
Differential Equations of Second Order,} Springer-Verlag, Berlin,
1983, MR0737190, Zbl 0562.35001.

\bibitem{gross} M.~Gross, ``Topological
Mirror Symmetry,'' Invent. Math. {\bf 144} (2001) 75--137,
MR1821145, Zbl pre01655606; and ``Special Lagrangian Fibrations I:
Topology,'' in {\sl Winter School on Mirror Symmetry, Vector
Bundles and Lagrangian Submanifolds,} C. Vafa and S.-T. Yau, eds.,
AMS/International Press (2001) 65--93, MR1876066, Zbl 0964.14033.


\bibitem{gs} M. Gross and B. Siebert,
``Affine Manifolds, Log Structures, and Mirror Symmetry,'' Turkish
Journal of Mathematics {\bf 27} (2003) 33-60; ``Mirror Symmetry
via Logarithmic Degeneration Data I,'' math.AG/0309070, MR1975331,
Zbl 1063.14048.

\bibitem{GW} M. Gross and P. M. H. Wilson,
``Large Complex Structure Limits of $K3$ Surfaces,'' J. Diff.
Geom. {\bf 55} (2000) 475-546. math.DG/0008018, MR1863732, Zbl
1027.32021.

\bibitem{hz02} C. Haase and I. Zharkov,
``Integral Affine Structures on Spheres and Torus Fibrations of
Calabi-Yau Toric Hypersurfaces I," math.AG/0205321.

\bibitem{hz03} C. Haase and I. Zharkov, ``Integral Affine
Structures on Spheres and Torus Fibrations of
Calabi-Yau Toric Hypersurfaces II," math.AG/0301222.

\bibitem{hit}
N. Hitchin, ``The Moduli Space of Special Lagrangian
Submanifolds,'' dedicated to Ennio De Giorgi, Ann. Scuola Norm.
Sup. Pisa Cl. Sci. {\bf 25}  (1997) 503--515, MR1655530, Zbl
1015.32022.

\bibitem{hitchin}  N. Hitchin, ``Lectures on Special Lagrangian
Submanifolds,'' in {\sl Winter School on Mirror Symmetry, Vector
Bundles and Lagrangian Submanifolds (Cambridge, MA, 1999)},
 AMS/IP Stud.\ Adv.\ Math.\ {\bf 23}, Amer.\ Math.\ Soc.\
(2001) 151--182, math.DG/9907034, MR1876068, Zbl pre01724275.

\bibitem{KS} M. Kontsevich and Y. Soibelman, ``Homological Mirror
Symmetry and Torus Fibrations,'' in {\sl Symplectic Geometry and
Mirror Symmetry,} World Scientific (2001) 203--263, MR1882331, Zbl
pre01787195.

\bibitem{lang85} S. Lang, {\sl Differential Manifolds,}
Springer-Verlag, New York, 1985, MR0772023, Zbl 0551.58001.

\bibitem{loftin02} J. Loftin, ``Affine Spheres and
K\"ahler-Einstein Metrics," Math. Res. Lett. \textbf{9} (2002)
425--432, MR1928863, Zbl 1033.53039.

\bibitem{loftin03} J. Loftin, ``The Compactification of the Moduli
Space of $\mathbb{RP}^2$ Surfaces, I," J. Diff. Geom. {\bf 68}
(2004) 223--276, math.DG/0311052, MR2144248.

\bibitem{loftin04} J. Loftin, ``Singular Semi-Flat Calabi-Yau
Metrics on $S^2$," Comm. Anal. Geom. {\bf 13} (2005) 333--361,
math.DG/0403218.

\bibitem{mcintosh} I. McIntosh, ``Special Lagrangian Cones in $\co^3$
and Primitive Harmonic Maps," J. London Math. Soc. (2) \textbf{67}
(2003) 769--789, MR1967705, Zbl pre02072641.

\bibitem{nomizu-sasaki} K. Nomizu and T. Sasaki, {\sl Affine Differential
Geometry,} Cambridge University Press, Cambridge, 1994, MR1311248,
Zbl 0834.53002.

\bibitem{o-v}  H. Ooguri and C. Vafa, ``Summing up Dirichlet Instantons,'' Phys.
Rev. Lett. \textbf{77} (1996) 3296--3298, MR1411842, Zbl
0944.81528.

\bibitem{wd-ruan} W.D. Ruan, ``Lagrangian Torus Fibration of Quintic Hypersurfaces. I. Fermat Quintic
Case," in {\sl Winter School on Mirror Symmetry, Vector Bundles
and Lagrangian Submanifolds (Cambridge, MA, 1999)}, AMS/IP (2001)
297--332, MR1876075, Zbl pre01724282; ``Lagrangian Torus Fibration
of Quintic Calabi-Yau Hypersurfaces. II. Technical Results on
Gradient Flow Construction," J. Symplectic Geom. \textbf{1} (2002)
435--521, MR1959057; ``Lagrangian Torus Fibration of Quintic
Calabi-Yau Hypersurfaces. III. Symplectic Topological SYZ Mirror
Construction for General Quintics," J. Diff. Geom. \textbf{63}
(2003) 171--229, MR2015547, Zbl pre02171927.

\bibitem{ruuska} V. Ruuska, ``Riemannian Polarizations," Ann. Acad.
Sci. Fenn. Math. Diss., \textbf{106} (1996), 38 pp., MR1413839,
Zbl 0862.53028.

\bibitem{simon-wang}
U.~Simon and C.-P. Wang.
\newblock Local theory of affine 2-spheres.
\newblock In {\em Differential Geometry: Riemannian geometry (Los Angeles, CA,
  1990)}, volume 54-3 of {\em Proceedings of Symposia in Pure Mathematics},
  Amer.\ Math.\ Soc.\ (1993) 585--598, MR1216648, Zbl 0804.53013.


\bibitem{syz}  A. Strominger, S.-T. Yau, and E. Zaslow, ``Mirror Symmetry is
T-Duality,'' Nuclear Physics \textbf{B479} (1996) 243-259,
MR1429831, Zbl 0896.14024.

\bibitem{tzitz} G. Tzitz\'eica, ``Sur une nouvelle classe de
surfaces," Rend. Circ. Mat. Palermo \textbf{25} (1908) 180--187,
JFM 39.0685.05; \textbf{28} (1909) 210--216, JFM 40.0668.04.

\bibitem{vafa}  C. Vafa, ``Extending Mirror Conjecture to Calabi-Yau with
Bundles,'' hep-th/9804131.

\bibitem{zh} I. Zharkov, ``Limiting Behavior of Local Calabi-Yau
Metrics," Adv. Theor. Math. Phys. \textbf{8} (2004) 395--420,
math.DG/0304116, MR2105186, Zbl pre02097019.

\end{thebibliography}

\vskip 0.1in

{\scriptsize John Loftin, Department of Mathematics, Rutgers
University, Newark, NJ 07102. (LOFTIN@ANDROMEDA.RUTGERS.EDU)}

{\scriptsize Shing-Tung Yau, Department of Mathematics, Harvard
University, Cambridge, MA  01238. (YAU@MATH.HARVARD.EDU) }

{\scriptsize Eric Zaslow, Department of Mathematics, Northwestern
University, Evanston, IL  60208.
(ZASLOW@MATH.NORTHWESTERN.EDU)}

\end{document}


\begin{center}
{\large \bf ERRATUM TO AFFINE MANIFOLDS, SYZ GEOMETRY AND THE ``Y"
VERTEX}\\[.1in]
JOHN LOFTIN, SHING-TUNG YAU, AND ERIC ZASLOW
\end{center}

\section{Main result}

The purpose of this erratum is to correct an error in the proof of
the main result of \cite{lyz05}. A semi-flat Calabi-Yau structure on
a smooth manifold $M$ consists of an affine flat structure, together
with an affine K\"ahler metric whose potential satisfies $\det
\Phi_{ij}=1$ in the local affine coordinates.  The main result of
\cite{lyz05} is
\begin{mr} There exist many
nontrivial semi-flat Calabi-Yau structures on the complement of a
trivalent vertex of a graph inside a ball in $\re^3$.
\end{mr}

We present two separate constructions of the Main Result below: The
first, presented in Section \ref{new-proof} below, uses hyperbolic
affine sphere structures constructed in \cite{loftin02c}.  The
second approach, proved in Section \ref{fix-thm2} below, is to use
elliptic affine spheres constructed by

\vspace*{.1in}

\noindent
 {\bf Theorem $\mathbf {2'}$.} {\it Let $U$ be a nonzero holomorphic
cubic differential on $\cp^1$ with exactly 3 poles of order two. Let
$M$ be $\cp^1$ minus the pole set of $U$. At each pole of $U$, let
$w$ be a local coordinate so that $U=\frac1{w^2}\,dw^3$ in a
neighborhood of the pole.  Let there be a conformal background
metric on $M$ which is equal to
$$\frac{|\log |w|^2|} {|w|}\, |dw|^2$$ on a neighborhood of each
pole. Then there is a $\delta>0$ so that for all
$\epsilon\in(0,\delta)$, there is a smooth bounded function $\eta$
on $M$ satisfying \begin{equation} \label{epsilon-U}\Delta \eta +
4\|\epsilon U \|^2 e^{-2\eta} + 2 e^\eta - 2\kappa = 0.
\end{equation}
(Here $\Delta$ is the Laplace operator of the background metric,
$\|\cdot \|$ is the induced norm on cubic differentials, and
$\kappa$ is the Gauss curvature.) The bound $\delta$ only depends on
$\int_M|U|^{\frac23}$.}

\vspace*{.1in}

This replaces Theorem 2 of \cite{lyz05}, which has a gap in its
proof.  The only essential difference between Theorem 2 of
\cite{lyz05} and Theorem $2'$ is that in Theorem 2, we claimed an
explicit value for the bound $\delta$, while now in Theorem $2'$, we
must assume $\delta$ is a small positive number without explicit
estimates.  We stress that this small strengthening of the
hypothesis does not affect the Main Result.

\section{Erratum}
The proof of Theorem 2 of \cite{lyz05} is incomplete. There is a
sign error in the proof of Part 2 of Proposition 8 which we have not
yet been able to overcome. We provide a separate proof of Theorem
$2'$, which is only slightly weaker, in Section \ref{fix-thm2}
below. We also give another alternate proof of the Main Result in
Section \ref{new-proof} below.

Theorem 2 was used to produce global solutions of the equation
governing two-dimensional elliptic affine spheres,
 \begin{equation} \label{eas} \Delta u + 4
\|U\|^2 e^{-2u} +2e^u - 2\kappa = 0,
\end{equation} on $\cp^1$ for $U$ a
nonzero cubic differential with three poles of order 2, and $\Delta$
and $\kappa$ the Laplacian and curvature of a background metric.
Theorem 2 claimed to find a solution of (\ref{eas}) for each $U$
which is small in the sense that
 \begin{equation} \label{sharp-bound} \int_{\cp^1} |U|^{\frac23}
< \frac{\pi\sqrt[3]2}{3}.\end{equation} Below in the proof of
Theorem $2'$, we produce solutions to (\ref{eas}) for all such $U$
satisfying
$$\int_{\cp^1} |U|^{\frac23}  < \delta $$ for a small $\delta>0$.
Then a result of Baues-Cort\'es \cite{baues-cortes03} relates
two-dimensional elliptic affine spheres to three-dimensional
semi-flat Calabi-Yau metrics, and the three poles of $U$ form the
rays of the ``Y" vertex. (The original bound $\frac{\sqrt[3]2}{3\pi}
$ is sharp if it is true \cite{lyz05}.)

The examples in Sections 2,4,5 of \cite{lyz05} and the construction
of local radially symmetric solutions to equation (\ref{eas}) near a
pole of $U$ in Section 6.1 of \cite{lyz05} are unaffected by Theorem
2.


\section{Hyperbolic affine spheres and metrics near the ``Y" vertex} \label{new-proof}
Our first alternative proof of the Main Result is based on the
following analogue of the result of Baues-Cort\'es:
\begin{prop} \label{hyp-as-cone}
Assume $n$ is even. Let $H$ be a hyperbolic affine sphere in
$\re^{n+1}$ centered at the origin with affine mean curvature $-1$.
Replace $H$ by an open subset of itself if necessary so that each
ray through the origin hits $H$ only once. Let $$ \mathcal C' =
\bigcup_{0<r<1} r H$$ be the union of all line segments from the
origin to $H$. Then the function $\Phi\!:\mathcal C' \to \re$ given
by
$$\Phi = \int(r^{n +1} - 1) ^{\frac 1{n+1}} \,dr$$ is convex
and solves the real Monge-Amp\`ere equation $\det \Phi_{ij} =
1$ on $\mathcal C'\subset \re^{n+1}$.
\end{prop}
The proof of this proposition is the same as our proof of
Baues-Cort\'es's result (which we label as Theorem 1 in
\cite{lyz05}).

Proposition \ref{hyp-as-cone} then reduces the problem of finding a
semi-flat Calabi-Yau structure on a neighborhood of the ``Y" vertex
of a graph in $\re^3$ to finding a hyperbolic affine sphere
structure on $S^2$ minus 3 singular points.   The next proposition,
which follows from \cite{loftin02c}, provides many such examples.

\begin{prop} \label{exist-has}
Let $U$ be a meromorphic cubic differential on $\cp^1$ with
three poles, each of order at most 3.  Then there is a
background metric $h$ and a solution $u$ to
\begin{equation} \Delta u + 4 \|U\|^2 e^{-2u} -2e^u - 2\kappa =
0
\end{equation}
so that $e^uh$ is the metric of a hyperbolic affine sphere.
\end{prop}

This proposition, together with a result of Wang \cite{wang91},
exhibits many hyperbolic affine sphere structures on $S^2$ minus the
pole set of $U$. Thus we find, using Proposition \ref{hyp-as-cone},
many nontrivial examples of semi-flat Calabi-Yau metrics on a
neighborhood of the ``Y" vertex of a graph in $\re^3$. Moreover, the
holonomy type around each puncture on $S^2$ (but not the global
holonomy representation) is determined in \cite{loftin02c}.

\section{Construction of elliptic affine spheres} \label{fix-thm2}
This section is devoted to the proof of Theorem $2'$, using the
barrier methods of \cite{loftin04}. First of all, we note
\begin{lem}
A function $\eta$ satisfies (\ref{epsilon-U}) if and only if
$\nu=\eta-\log \epsilon$ satisfies
\begin{equation}
\label{perturb-equation}
 \Delta \nu + 4\|U\|^2 e^{-2\nu} + 2\epsilon e^\nu - 2\kappa = 0.
\end{equation}
\end{lem}

It will be convenient to pass to a branched cover of $\cp^1$ under
which the pullback of $U$ has 6 poles of order one.  For example,
consider the map $$\Pi \!: Z \mapsto Z^2 + Z^{-2},$$ for $Z$ an
inhomogeneous coordinate on $\cp^1$.  $\Pi$ is then a branched cover
from $\cp^1\to\cp^1$ with branch points over $Z= 2,-2,\infty$.  If
we allow the (order-two) poles of $U$ to be at $2,-2,\infty$, then
$V=\Pi^*U$ has six poles of order one at $0,\infty,1,-1,i,-i$.  Let
$G$ be the group of deck transformations of $\Pi$. Note $|G|=4$.

This discussion shows that Theorem $2'$ is equivalent to the
existence of a $G$-invariant bounded solution $\mu$ of
\begin{equation} \label{our-eq}
L_{h,\epsilon}(\mu)=  \Delta_h \mu + 4\|V\|_h^2 e^{-2\mu} +
2\epsilon e^\mu - 2\kappa_h = 0
\end{equation}
for $\epsilon<\delta$ and for a $G$-invariant background metric $h$
on $N=\Pi^{-1}M $ which is equal to $|\log|z|^2| \, |dz|^2$ near the
poles of $V$. (Here the coordinate $z$ is chosen so that $V=\frac1z
dz^3$ near $z=0$.)

\begin{prop}
There is a $\delta>0$ so that there are smooth bounded $G$-invariant
functions $s, S$ on $N$ satisfying
\begin{eqnarray}
 s &\le & S, \\
 L_{h,\epsilon}(s) & \ge & 0, \\
 L_{h,\epsilon}(S) &\le &0
\end{eqnarray}
for all $\epsilon\in(0,\delta)$.
\end{prop}

\begin{proof}
 These barriers follow from the same arguments as in Section 3.1 of
 \cite{loftin04}. Recall $z$ is a local coordinate near each pole of
 $V$ so that $V = \frac1z \,dz^3$ near the pole $z=0$, and that the
background metric is equal to $|\log |z|^2| \,|dz|^2$ near $z=0$.
All the explicit formulas for barriers below are invariant under the
action of $G$:  a nontrivial element of $G$ acts at a fixed point
$z=0$ by $z\mapsto -z$, and the formulas below only depend on $|z|$.

Let $\alpha\in (-1,0)$, and let $f$ be a $G$-invariant positive
function which is equal to $|\log |z| |^\alpha$ in a neighborhood of
each pole of $V$ and is equal to a constant outside a neighborhood
of the poles.  Then by the calculations in \cite{loftin04}, for
$\beta\ll0$, $s = \beta f$ satisfies
$$ L_{h,\epsilon}(s) \ge L_{h,0}(s) \ge 0.$$
(This computation uses the fact that $V$ has no zeros, although the
case in which $V$ has zeros can be handled, too \cite{loftin04}.)

On the other hand, to produce $S$, let $k$ be a smooth $G$-invariant
metric on $\cp^1$ which is equal to $|dz|^2$ near each pole of $V$.
Let $\tilde\kappa$ be a $G$-invariant positive function on $N$ which
is equal to
$$\frac 1 { 2 |z|^2 (\log |z|)^2}$$ in a neighborhood of each pole
and so that $$\int_N (\kappa_k - \tilde \kappa) dV_k = 0$$ for
$dV_k$ the volume form of $k$. Also, let $\xi$ be a $G$-invariant
function equal to $\log|\log|z|^2|$ near each pole $z=0$ and smooth
away from the poles. Apply the Green's function of $\Delta_k$ to
produce a $G$-invariant function $\tilde g$ so that
$$\Delta _k \tilde g = 2\kappa_k - 2\tilde \kappa -\Delta_k\xi.$$
Since $\Delta_k (\log|\log|z|^2|) = -\frac1{|z|^2(\log|z|)^2} =
-2\tilde\kappa$ in the sense of distributions near $z=0$, $\tilde g$
is harmonic near $z=0$ and is thus smooth. By construction, $g =
\tilde g + \xi$ satisfies
$$g = \log |\log |z|^2| + O(1)$$ near each pole of $V$. (We may also
produce $g$ by applying the Green's function to $2\kappa_k - 2\tilde
\kappa$ as in \cite{loftin04}.) Compute for a constant $c$,
$$ L_{k,0}( g+c) = 4 e^{-2g-2c} \|V\|^2_k - 2\tilde \kappa < 0$$ on
$N$ for $c\gg0$, since $e^{-2g}\|V\|_k^2$ and $\tilde\kappa$ have
the same order of growth as $z\to0$. If we relate the metrics $h$
and $k$ by $k = e^v h$, then compute for $S=g+c+v$
$$L_{h,0} (S) = e^v L_{k,0}(g+c) < 0.$$
The asymptotics of $g$ and the definition of $v$ ensure that $S$ is
bounded. Moreover,
$$L_{h,0}(S) \le - \frac C {|z|^2 |\log |z| |^3}$$ for a constant
$C>0$ near each pole of $V$.  On the other hand, $e^S$ is bounded,
and so on all of $N$ $$L_{h,\epsilon}(S) = L_{h,0}(S) + 2\epsilon
e^S < 0$$ for all positive $\epsilon$ which are small enough.
\end{proof}

Now we use the method of sub- and super-solutions to produce a
$G$-invariant solution to $L_{h,\epsilon} (\mu)=0$ on $N$. Exhaust
$N$ by compact submanifolds $N_n$ with boundary $$ N = \bigcup_n
N_n, \qquad N_n \subset \mbox{int} \,N_{n+1}.$$ Assume each $N_n$ is
$G$-invariant. We proceed as in e.g.\ Schoen-Yau \cite{schoen-yau94}
Proposition V.1.1 (which is essentially the same result on compact
manifolds without boundary)
\begin{prop}
For $\epsilon$ sufficiently small, there exists a $G$-invariant
solution $\mu_n$ to (\ref{our-eq}) on $N_n$ which satisfies $s\le
\mu_n\le S$ on the interior of $N_n$.
\end{prop}
\begin{proof}
Write $L_{h,\epsilon}(\mu) = \Delta_h \mu + f(x,\mu)$ for $f(x,t) =
4\|V\|_h^2 e^{-2t} + 2\epsilon e^t - 2\kappa_h$ and $x\in N_n$.   By
compactness, there is a constant $\gamma>0$ so that $F(x,t) = \gamma
t + f(x,t)$ is increasing for $t\in [-A,A]$ for a constant $A=A_n$
satisfying $-A \le s \le S \le A$ on $N_n$ (for now, we suppress the
dependence on $n$). Define an operator $H$ by $H \mu = -\Delta_h \mu
+ \gamma\mu$, and define inductively $s_k$ and $S_k$
$$\begin{array}{c}
 s_0 = s, \qquad H(s_k) = F(x,s_{k-1}), \quad s_k|_{\partial N_n} =
 S. \\
 S_0 = S, \qquad H(S_k) = F(x,S_{k-1}), \quad S_k|_{\partial N_n } =
 S.
\end{array}$$
Note that the Dirichlet problems admit unique solutions to make
$s_k,S_k$ well-defined. The functions $s_k,S_k$ are $G$-invariant by
this uniqueness.  As in \cite{schoen-yau94}, we use the maximum
principle to show that
$$s=s_0\le s_1 \le s_2 \le \cdots \le s_k \le\cdots \le S_k \le
\cdots \le S_2 \le S_1 \le S_0 = S.$$ Then the monotone pointwise
limits $s_k \to s_\infty$ and $S_k\to S_\infty$ are each
$C^{2,\alpha}$ solutions to (\ref{our-eq}) on the interior of $N_n$
(use the standard elliptic theory as in \cite{schoen-yau94}). Let
$\mu_n = s_\infty$, which  is $G$-invariant, since it is the limit
of $G$-invariant functions.
\end{proof}

\begin{prop}
As $n\to\infty$, there is a subsequence $n_i$ so that $\mu_{n_i} \to
\mu$ in $C^2_{\rm loc}$ on $N$.  $\mu$ is a $G$-invariant bounded
solution to (\ref{our-eq}).
\end{prop}

\begin{proof}
Choose $n$ and let $m\ge n$. Then $\mu_m$ are uniformly bounded
solutions to (\ref{our-eq}) on $N_n$, as $|\mu_m|\le A_n$ on $N_n$.
So there are $L^p$ estimates for $\mu_m$ and $\Delta_h \mu_m$ on
$N_n$, which lead to local $L^p_2$ estimates on $\mu_m$ on each open
subset of $\mathcal O_n \subset N_n$. Then we may use Sobolev
embedding to get $C^{0,\alpha}$ estimates on $\mu_m$ and $\Delta_h
\mu_m$:
$$\|\mu_m\|_{C^{0,\alpha}}, \|\Delta_h \mu_m \|_{C^{0,\alpha}} \le C_n
\mbox{ on } \mathcal O_n,$$ where we have shrunk $\mathcal O_n$
slightly. Then by the Schauder theory, there are uniform
$C^{2,\alpha}$ estimates for $\mu_m$ on $\mathcal O_n$ (which is
again shrunk a little). All these estimates are uniform in $m\ge n$
on each  $\mathcal O_n \subset N$.  We may still assume $\mathcal
O_n$ is still an exhaustion of $N$:
$$N = \bigcup_n \mathcal O_n, \qquad \mathcal O_n \Subset \mathcal
O_{n+1}.$$

 Then we may use the
Ascoli-Arzel\'a Theorem and the usual diagonalization argument to
produce a subsequence $\mu_{n_i} \to\mu$ converging in $C^{2}$ on
any compact subset of $N$, and so $\mu$ is a solution to
(\ref{our-eq}) bounded between $s$ and $S$. $\mu$ is $G$-invariant,
since it is a limit of functions which are all $G$-invariant.
\end{proof}

To complete the proof of Theorem $2'$, we simply note that all the
estimates only depend on $\|U\|$, and since the space of admissible
$U$ on $\cp^1$ has one complex dimension, $\delta$ depends only on
$\int_M |U|^{\frac23}$.

\begin{rem}
It is presumably possible to solve (\ref{epsilon-U}) for a
nontrivial cubic differential $U$ on $\cp^1$ which admits $\le3$
poles of order 2 and an arbitrary number of poles of order 1.  It is
not possible to have more than 3 poles of order 2, since the
curvature $\kappa$ of the background metric contributes a point mass
at each pole of order 2 (while the corresponding background metric
for poles of order 1 makes no such contribution).  We have
$$\int_{\cp^1} 2\kappa \,dV= 8\pi - 2\pi \ell$$ for $\ell$ the number of
poles of order 2.  Integrating equation (\ref{epsilon-U}) then shows
that there are no nontrivial solutions of (\ref{epsilon-U}) if
$\ell>3$. (In the present case, the essential fact is that the
orbifold $\cp^1/G$ admits a positive-curvature orbifold metric.)
\end{rem}

\bibliographystyle{abbrv}